\documentclass[12pt]{amsart}
\usepackage[utf8]{inputenc}
\usepackage[english]{babel}
\usepackage{subcaption}
\usepackage{enumerate}
\usepackage{amsfonts}
\usepackage{comment}
\usepackage{amsmath,amssymb}
\usepackage[makeroom]{cancel}
\usepackage{stmaryrd}
\usepackage{tikz-cd} 
\usepackage{enumerate}
\usepackage{enumitem}
\usepackage{graphicx}
\usepackage{amsthm}
\usepackage{wrapfig}
\usepackage{hyperref}
\usepackage{chngcntr}
\counterwithin{figure}{section}
\newtheorem{theorem}{Theorem}
\newtheorem*{theorem*}{Theorem}
\newtheorem{prop}[theorem]{Proposition}

\newtheorem{lemm}[theorem]{Lemma}
\theoremstyle{definition}
\newtheorem{defi/}[theorem]{Definition}
\newenvironment{defi}
  {%
   \pushQED{\qed}\begin{defi/}}
  {\popQED\end{defi/}}
\newtheorem*{defi*}{Definition}
\newtheorem*{prop*}{Proposition}
\newtheorem*{lemm*}{Lemma}

\newtheorem{prdef}[theorem]{Proposition-Definition}

\newtheorem{cor-definition}[theorem]{Corollary-Definition}
\newtheorem{rema}[theorem]{Remark}

\newtheorem*{conv*}{Convention}

\numberwithin{theorem}{section}

\newcommand{\inte}[1]{\overset{\circ}{#1}}

\newcommand{\quotient}[2]{{\raisebox{.2em}{$#1$}\left/\raisebox{-.2em}{$#2$}\right.}}

\numberwithin{equation}{section}
\begin{document}
\title{Dehn-Fried surgeries on non-transitive expansive flows} 
\author[I. Iakovoglou]{Ioannis Iakovoglou}
\address{Université Sorbonne Paris Nord \\ \vspace{0.1cm} UMR 7539 du
CNRS, 99 avenue J.B. Clément, 
93430 Villetaneuse, 
FRANCE.}
\email{e.mail: ioannis.iakovoglou@ens-lyon.fr}
\begin{abstract}
 In this paper, we prove that for any flow $\Psi$ obtained by a Dehn-Fried surgery on an expansive (or equivalently topological pseudo-Anosov) flow $\Phi$ in dimension 3, $\Psi$ is expansive if and only if its stable and unstable foliations do not contain $1$-prong singularities. 
\end{abstract}
\maketitle

\section{Introduction}
    Expansive flows in dimension 3 form a large family of dynamical systems with several hyperbolic properties. More specifically, thanks to a classical result by Inaba,  Matsumoto and Oka, any expansive flow in dimension 3 preserves a pair of transverse singular foliations, the stable and unstable foliations, each singularity of which is a circle $k$-prong singularity with $k\geq 2$. 

    Examples of expansive flows in dimension 3 include : 1) suspensions of pseudo-Anosov homeomorphisms on closed surfaces 2) Anosov flows 3) certain flows that are transverse to foliations in closed atoroidal manifolds (see \cite{Cal1}, \cite{Cal2}, \cite{Cal3}, \cite{Fe1}, \cite{Mo}) 4) expansive flows obtained by Dehn-Fried surgeries along periodic orbits of transitive pseudo-Anosov flows (see \cite{Fried}, \cite{Mariothese}). 
    
    The Dehn-Fried surgery constitutes nowadays one of the main tools for constructing new transitive pseudo-Anosov flows in dimension 3. It consists in performing a Dehn surgery on the tubular neighborhood of some periodic orbit of a transitive pseudo-Anosov flow in order to create a new transitive pseudo-Anosov flow supported by a new manifold (see Definition \ref{d.defisurgery} for more details). Thanks to the existence of Birkhoff sections for every transitive pseudo-Anosov flow (see \cite{Brunella}), it is not difficult to show that Dehn-Fried surgeries preserve the class of transitive expansive flows. In this paper, we will prove that this result also holds, in most cases, for non-transitive expansive flows in dimension 3. 

    More specifically, consider an expansive flow $\Phi$ in dimension 3 and $F_{\Phi}^s,F_{\Phi}^u$ its stable and unstable foliations. Let $\Psi$ be any flow obtained by a Dehn-Fried surgery along a periodic orbit of $\Phi$, whose tubular neighborhoods are homeomorphic to solid tori. The previous surgery carries $F_{\Phi}^s,F_{\Phi}^u$ into a pair of (possibly singular) transverse foliations $F_{\Psi}^s,F_{\Psi}^u$ (see Proposition \ref{p.foliationsaftersurgery}). 
    \begin{theorem}\label{t.maintheorem}
         The flow $\Psi$ is expansive if and only if $F_{\Psi}^s,F_{\Psi}^u$ do not admit a circle one prong singularity. 
    \end{theorem}
Even though nowadays the previous result is considered by many as folklore, the author of this paper considers the proof of the above theorem to be an important omission in the theory of pseudo-Anosov flows in dimension 3, as to our knowledge, every known example of non-transitive pseudo-Anosov (and not Anosov) flow in dimension 3, with the exception of the flows constructed in \cite{BarbotFenley}, is constructed via a Dehn-Fried surgery. 

\textit{Acknowledgments}. The completion of this paper was made possible thanks to the grant ``ERC Emergence" that financed my postdoctoral fellowship in Sorbonne Université under the supervision of Pierre Berger. 
\section{Preliminaries}\label{s.prelim}
\subsection{Expansive flows and singular foliations}\label{ss.singfol}
\begin{defi}\label{d.expansiveflow}
    Let $M$ be a $C^0$ closed $3$-manifold endowed with a distance $d$ (compatible with its topology). Consider $\Phi=(\Phi^t)_{t\in \mathbb{R}}$ a non-singular $C^0$-flow on $M$ (i.e. a one parameter family of homeomorphisms verifying $\Phi^{t+s}=\Phi^t\circ \Phi^s$ and $\Phi^0=id_M$). We will say that $\Phi$ is \emph{expansive} if there exist $\epsilon, \eta>0$ such that for any $x,y\in M$ and $h:\mathbb{R}\rightarrow \mathbb{R}$ an increasing  homeomorphism with $h(0)=0$ $$ \forall t\in \mathbb{R}~ d(\Phi^t(x),\Phi^{h(t)}(y))<\eta \implies \exists |t_0|<\epsilon ~~ y=\Phi^{t_0}(x)$$
\end{defi}
Note that, thanks to the compactness of $M$, the expansive character of a flow does not depend on our choice of metric. 

Consider $\mathcal{F}^s_2$ (resp. $\mathcal{F}^u_2$) the foliation by horizontal (resp. vertical) lines on $\mathbb{C}$, $\mathcal{D}_2$ the euclidean square $\{z\in\mathbb{C}||\text{Re}(z)|<1 \text{, } |{Im}(z)|<1\}$ and $\pi_p(z)=z^p$ for every $p\in \mathbb{N}^*$. Let   $\mathcal{D}_1= \pi_2(\mathcal{D}_2)$ and $\mathcal{D}_p= \pi_p^{-1}(\mathcal{D}_1)$ for any $p\geq 3$. 

The image of $\mathcal{F}^s_2$ (resp. $\mathcal{F}^u_2$) by $\pi_2$ defines a singular foliation $\mathcal{F}^s_1$ (resp. $\mathcal{F}^s_1$). Similarly, for every $p\geq 3$ the pre-image of $\mathcal{F}^s_1$ (resp. $\mathcal{F}^u_1$) by $\pi_p$ defines a singular foliation, say $\mathcal{F}^s_p$ (resp. $\mathcal{F}^u_p$). We will say that $\mathcal{F}^s_1$ (or  $\mathcal{F}^u_1$) has a \emph{$1$-prong singularity at $0$} and that $\mathcal{F}^s_p$ (or  $\mathcal{F}^u_p$) has a \emph{$p$-prong singularity at $0$}.

\begin{defi}\label{d.folisingular}
   Consider $M$ a closed $3$-manifold. We will say that $\mathcal{F}$ is \emph{a singular codimension one foliation of $M$} if it is a decomposition of $M$ into \emph{regular} and \emph{singular leaves} satisfying the following properties: 
    \begin{enumerate}
    \item there exists a finite number (possibly zero) of singular leaves $L_1,...,L_n$ 
    \item for every $i\in \llbracket 1,n\rrbracket$ there exists  $C_i\subset L_i$ an embedded circle in $M$ such that $$\text{regular leaves}\cup \{L_1-C_1\}\cup...\cup \{L_n-C_n\}$$ defines on $M-\overset{n}{\underset{i=1}{\cup}}C_i$ a $C^0$ codimension one  foliation. 
        \item for every $x\in C_i$ there exist $U_x$ a neighborhood of $x$ in $M$, $p\in \mathbb{N}\setminus\{0,2\}$ and $h:U_x\rightarrow \mathcal{D}_p\times [0,1]$ a homeomorphism verifying $h(x)=(0,\frac{1}{2})$ and $h(\mathcal{F}\cap U_x)=\mathcal{F}_p^s\times [0,1]$
    \end{enumerate}
    We will call $C_i$ a \emph{circle $p$-prong singularity of $\mathcal{F}$} and we will denote the set of circle prong singularities of $\mathcal{F}$ as $\text{Sing}(\mathcal{F})$. 
\end{defi}
\begin{defi}\label{d.transversedim3}
    Take $\mathcal{F},\mathcal{G}$ two singular codimension one foliations on $M$. We will say that $\mathcal{F}$ and $\mathcal{G}$ are \emph{transverse} if:
    \begin{enumerate}
        \item $\text{Sing}(\mathcal{F})=\text{Sing}(\mathcal{G})$
        \item for every point $x\in M$ there exists $U_x$ a neighborhood of $x$ in $M$, $p\geq 1$ and $h:U_x\rightarrow \mathcal{D}_k\times [0,1]$ a homeomorphism such that 
        \begin{enumerate}
            \item $h(x)=(0,\frac{1}{2})$
            \item $h(\mathcal{F}\cap U_x)=\mathcal{F}_p^s\times [0,1]$ and $h(\mathcal{G}\cap U_x)=\mathcal{F}_p^u\times [0,1]$
        \end{enumerate} 
 \end{enumerate}
\end{defi}
\begin{defi}[Topological pseudo-Anosov flow]\label{d.pseudoanosovflow}
  Let $M$ be a closed $3$-manifold, $d$ a distance on $M$ (compatible with its topology) and $\Phi=(\Phi^t)_{t\in\mathbb{R}}$ a $C^0$-flow on $M$. We will say that $\Phi$ is an \emph{almost pseudo-Anosov} flow if: 
  \begin{enumerate}
      \item $\Phi^t$ is non-singular 
      \item $\Phi^t$ preserves a pair of transverse singular codimension one foliations $F^s$, $F^u$ with no circle 1-prong singularities, called respectively the \emph{stable} and \emph{unstable foliations}
      \item for every $x\in M,y\in F^s(x)$ (resp. $y\in F^u(x)$) there exists an increasing homeomorphism $h:\mathbb{R}\rightarrow \mathbb{R}$ such that $$d(\Phi^t(x),\Phi^{h(t)}(y))\underset{t\rightarrow +\infty}{\longrightarrow} 0~~ \big(\text{resp. }d(\Phi^t(x),\Phi^{h(t)}(y))\underset{t\rightarrow -\infty}{\longrightarrow}0\big) $$
   \end{enumerate}   
   If furthermore $\Phi$ satisfies the following condition, we will say that $\Phi$ is a \emph{topological pseudo-Anosov} flow: 
   \begin{enumerate}\setcounter{enumi}{3}
      \item there exist $\eta, \epsilon>0$ such that for any $x\in M$, $y\in F^s(x)$, $z\in F^u(x)$ and any increasing  homeomorphism $h:\mathbb{R}\rightarrow \mathbb{R}$ with $h(0)=0$ $$ \forall t\geq 0~ d(\Phi^t(x),\Phi^{h(t)}(z))<\eta \implies \exists |t_0|<\epsilon ~~ z=\Phi^{t_0}(x)$$
      $$ \forall t\leq 0~ d(\Phi^t(x),\Phi^{h(t)}(y))<\eta \implies \exists |t_0|<\epsilon ~~y=\Phi^{t_0}(x)$$
  \end{enumerate} 
\end{defi}
Once again, thanks to the compactness of $M$, the property of being an almost or a topological pseudo-Anosov flow does not depend on our choice of metric $d$. By Theorem 1.5 of \cite{Inaba}\footnote{The proof of Theorem 1.5 in \cite{Inaba} requires the underlying manifold to be orientable. However, the orientability hypothesis is not used in the proof of the existence of the stable and unstable foliations.} and Lemma 2.7 of \cite{Oka} we have that: 
\begin{theorem}[Inaba, Matsumoto, Oka]\label{t.expansiveimpliespseudo}
  Let $M$ be a closed $3$-manifold and $\Phi$ an expansive flow on $M$. The flow $\Phi$ is a \emph{topological pseudo-Anosov flow}. 
\end{theorem}

Conversely, it is a classical result that 
\begin{theorem}\label{t.pseudoimpliesexpansive}
Any topological pseudo-Anosov flow on a 3-manifold is an expansive flow. 
\end{theorem}

\subsection{Polygons and local models of periodic orbits}\label{s.polygonsandmodels}
Define $ \pi_p, \mathcal{D}_p,\mathcal{F}^s_p, \mathcal{F}^u_p$ for every $p\geq 1$ as in Section \ref{ss.singfol}. 

\begin{defi}
    For any $p\geq 1$ and any $x\in \mathbb{R}^2$, the closure of any connected component of $\mathbb{R}^2-(\mathcal{F}^s_p(x)\cup \mathcal{F}^u_p(x))$ will be called a \emph{quadrant} of $x$ in $(\mathbb{R}^2,\mathcal{F}^s_p,\mathcal{F}^u_p)$. Moreover, the closure of any connected component of $\mathcal{F}^s_p(0)-\{0\}$ (resp. $\mathcal{F}^u_p(0)-\{0\}$) will be called a stable (resp. unstable) \emph{prong} of the origin. 
\end{defi} 

Denote by $\phi_2:\mathbb{R}^2\rightarrow \mathbb{R}^2$ the map $\phi_2(x,y)= (\frac{1}{2}x, 2 y)$. Notice that the foliations $\mathcal{F}^s_2, \mathcal{F}^u_2$ form the stable and unstable foliations of $\phi_2$. Furthermore, $\phi_2$ projects to a unique homeomorphism $\phi_1:\mathbb{R}^2\rightarrow \mathbb{R}^2$ preserving $\mathcal{F}^s_1$ and $\mathcal{F}^u_1$ and satisfying $\pi_2\circ \phi_2=\phi_1\circ\pi_2$. For any $p\geq 3$, denote by $\phi_p:\mathbb{R}^2\rightarrow \mathbb{R}^2$ the unique homeomorphism preserving each prong of the origin inside $\mathcal{F}^s_p$ and $\mathcal{F}^u_p$ and satisfying $\pi_p\circ \phi_p=\phi_1\circ\pi_p$.

We define for every $p\geq 1$ and $k\in \llbracket 0, p-1\rrbracket $ $$\phi_{pk}:=\phi_p\circ R_{k/p}=R_{k/p} \circ \phi_p$$ where $R_{\theta}:\mathbb{R}^2\rightarrow \mathbb{R}^2$ is the rotation around the origin of angle $2\pi\theta$. Following the terminology of Mosher (see \cite{Mo}), we will call $\phi_{pk}$ the \emph{local model for a
pseudo-hyperbolic fixed point with $p$ prongs and rotation $k$}. By suspending the previous local model, we will now define a local and smooth model for a $p$-prong circle singularity. Let $p\geq 1$ and $k\in \llbracket 0, p-1\rrbracket $. Consider the mapping torus
$$N_{pk}:= \frac{\mathbb{R}^2\times \mathbb{R}}{((x,y),t+1) \sim (\phi_{pk}(x,y),t) }$$
endowed with the constant speed vertical flow $\Phi_{pk}=(\Phi^t_{pk})_{t\in \mathbb{R}}$ given by the vector field $\frac{\partial}{\partial t}$. We will call $(N_{pk},\Phi_{pk})$ the \emph{local model for a pseudo-hyperbolic periodic orbit with $p$ prongs and rotation $k$}. We will denote by $\gamma_{pk}$ the projection of $\{(0,0)\}\times \mathbb{R}$ on $N_{pk}$ ($\gamma_{pk}$ is the unique periodic orbit of $\Phi_{pk}$) and by $F^s_{pk},F^u_{pk}$ the projection of $\mathcal{F}^s_p\times \mathbb{R}, \mathcal{F}^u_p\times \mathbb{R}$ on $N_{pk}$. We will call $F^s_{pk},F^u_{pk}$ the \emph{stable and unstable foliations of $\Phi_{pk}$}. 

\begin{figure}
 \begin{minipage}[ht]{0.4\textwidth}
    \centering 
     \vspace{0.8cm}
    \hspace{-0.5cm}
    \includegraphics[width=0.28\textwidth]{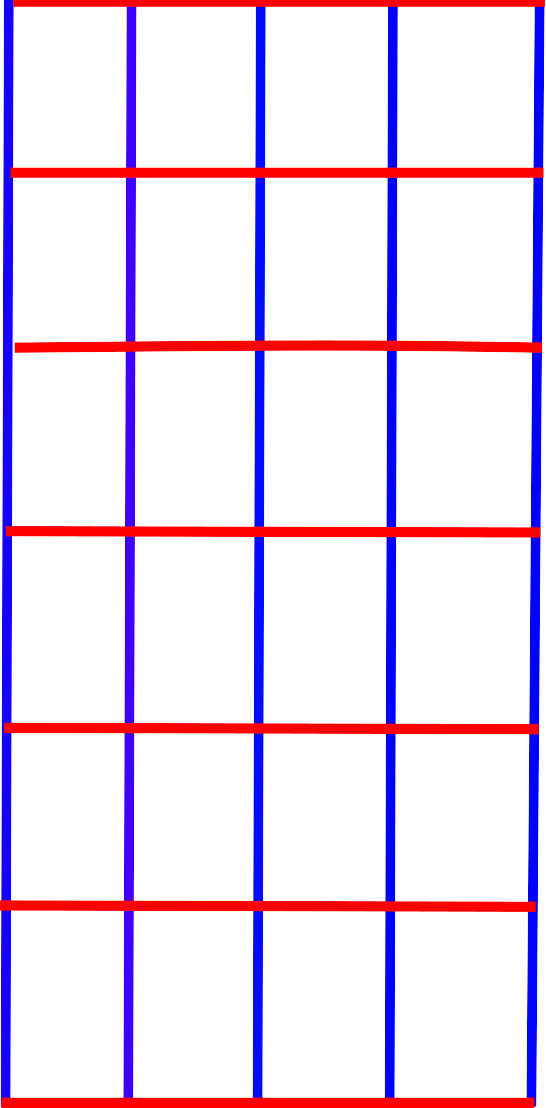}
  \hspace{-1cm}
    \caption*{\quad (a)}
    
  \end{minipage}
 \begin{minipage}[ht]{0.4\textwidth}
 \centering
 \vspace{0.9cm}
    \includegraphics[width=0.55\textwidth]{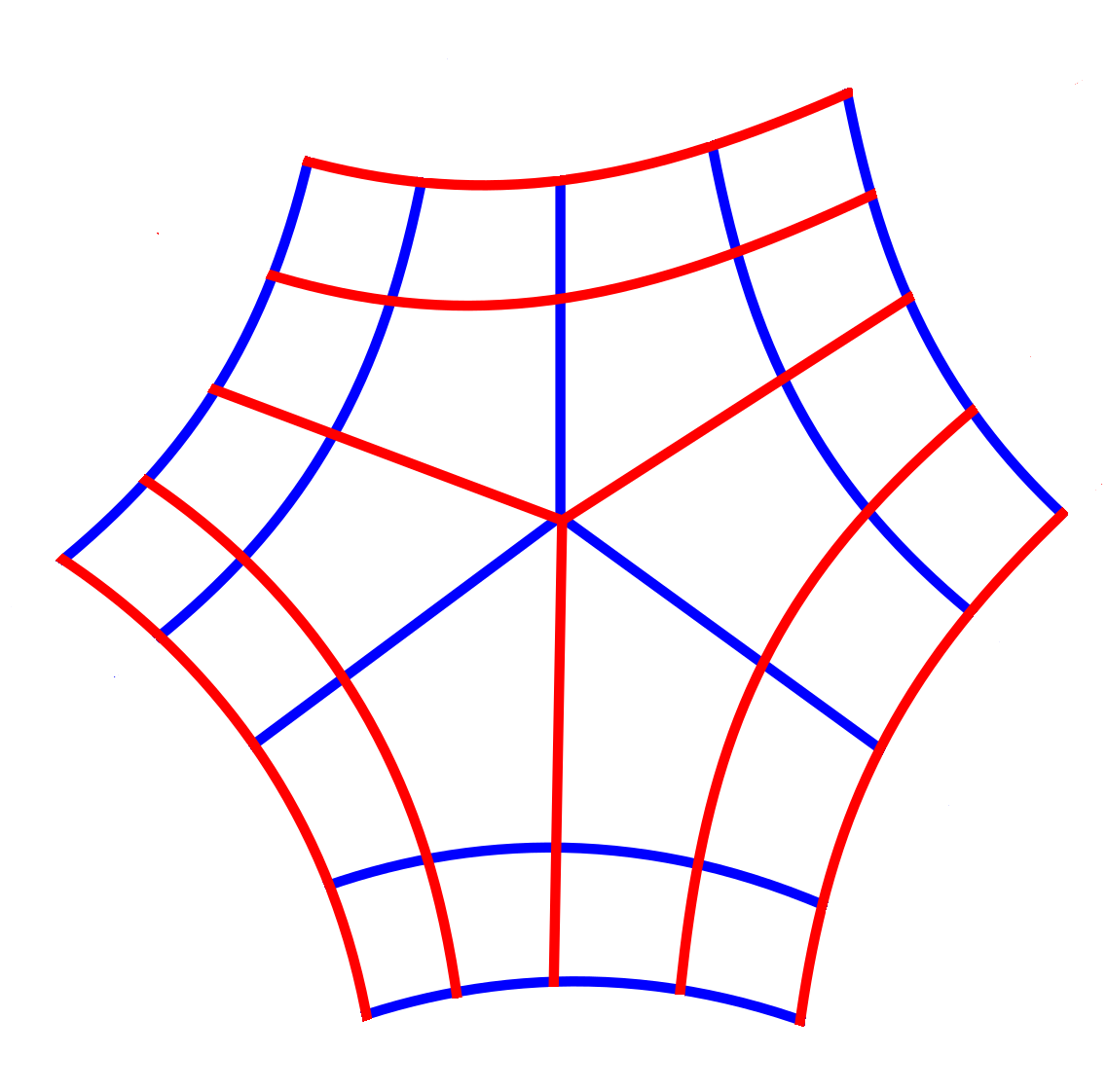}
    
    \caption*{(b)}
    
  \end{minipage}
\caption{(a) A rectangle (b) A standard hexagon}
\label{f.polygons}
  \end{figure}

Consider $\Phi$ an almost pseudo-Anosov flow on $M^3$, $F^s,F^u$ its stable and unstable foliations and $\gamma$ a periodic orbit of $\Phi$, whose tubular neighborhoods are homeomorphic to solid tori. Local models for pseudo-hyperbolic orbits describe the local behavior around every periodic orbit of $\Phi$ with  orientable tubular neighborhoods (or equivalently tubular neighborhoods that are homeomorphic to solid tori). This result relies essentially on the existence of standard polygons for $\Phi$: 
\begin{defi}\label{d.standardpolygonflow}
  A \emph{transverse standard polygon} $P$ in $M$ (for $\Phi$) is an embedded closed topological disk such that: 
 \begin{enumerate}
     \item $P$ is \emph{topologically transverse} to $\Phi$, i.e there exists $\epsilon>0$ such that for every $x\in P$ the segment $\underset{t\in(-\epsilon, \epsilon)}{\cup}\Phi^{t}(x)$ intersects $P$ once
     \item $P$ contains at most one point belonging to a circle prong singularity of $\Phi$
     \item when $\inte{P}:=P-\partial P$ intersects a circle $k$-prong singularity, $P$ is bounded by $k$ stable and $k$ unstable segments that alternate between them (see Figure \ref{f.polygons}). In this case, we call $P$ a \emph{transverse (standard) $2k$-gon}
     \item when $\inte{P}$ does not intersect a circle prong singularity, $P$ is bounded by $2$ stable and $2$ unstable segments that alternate between them. In this case we will call $P$ a \emph{transverse (standard) rectangle}.

\end{enumerate}   
\end{defi}
By Theorem 17A of \cite{Whitney} and the fact that $\Phi$ preserves a pair of transverse singular
foliations, it is easy to see that any $x\in M$ is contained in the interior of a transverse
standard polygon. Thanks to this fact, it is not hard to prove that:
\begin{prop}\label{p.aroundcircleprong}
    Consider $(N_{pk},\Phi_{pk})$ the local model for a pseudo-hyperbolic periodic orbit $\gamma_{pk}$ with $p\geq 2$ prongs and  rotation $k \in \llbracket 0, p-1 \rrbracket$. For every periodic orbit $\gamma$ of $\Phi$, whose tubular neighborhoods are solid tori, there exist $p\geq 2$, $k \in \llbracket 0, p-1 \rrbracket$, $U_\gamma$ a neighborhood of $\gamma$ homeomorphic to a solid torus, $V_{pk}$ a neighborhood of $\gamma_{pk}$ and a homeomorphism $H:U_{\gamma}\rightarrow V_{pk}$ such that 
    \begin{enumerate}
    \item $H$ defines an orbital equivalence between the restriction of $\Phi$ on $U_{\gamma}$ and the restriction of $\Phi_{pk}$ on $V_{pk}$ 
        \item for any $x\in U_{\gamma}$ the map $H$ takes the connected component of $F^s(x)\cap U_\gamma$ (resp. $F^u(x)\cap U_\gamma$) containing $x$ to the connected component of $F^s_{pk}(H(x))\cap V_{pr}$ (resp. $F^u_{pk}(H(x))\cap V_{pr}$) containing $H(x)$  
    \end{enumerate}
\end{prop}

\section{The definition of the Dehn-Fried surgery}
Following the original definition of Fried (see \cite{Fried}), any Dehn-Fried surgery on an expansive flow $\Phi$ consists of two distinct geometric/dynamical operations: 
\begin{enumerate}
    \item the \emph{blow-up operation}, which consists in exploding a periodic orbit $\gamma$ (with orientable tubular neighborhoods) to a torus. In more geometric terms, 
    if the underlying manifold of $\Phi$ and the curve $\gamma$ were smooth, then the blow-up operation consists in replacing every point in $\gamma$ by the unit normal bundle of $\gamma$ at that point
    \item the \emph{blow-down operation}, which consists in crushing the torus obtained after the previous explosion to a new circle
\end{enumerate}
In the following pages, we will define in detail each of the previous operations first for hyperbolic models and then for general expansive flows. 
\subsection{Blow-up operation for hyperbolic models}\label{s.blowupmodels}
Define $ \pi_p, \phi_p, \phi_{pk}, \mathcal{D}_p,\mathcal{F}^s_p, \mathcal{F}^u_p, N_{pk}$, $\gamma_{pk}, F^s_{pk}, F^u_{pk}, R_{\theta}$ for every $p\geq 1$, $k\in \llbracket 0, p-1\rrbracket $ and $\theta \in [0,1]$ as in Section \ref{s.prelim}. For every $p\in \mathbb{N}\setminus\{0,2\}$ the map $\phi_p$ defines a homeomorphism that is $C^{\infty}$ at every point of $\mathbb{R}^2$ except at the origin ($\phi_2$ is everywhere smooth). Despite its lack of differentiability at the origin, by construction $\phi_p$ acts on the set of half-lines based at the origin. It follows that after blowing up the origin in $\mathbb{R}^2$, we obtain a surface diffeomorphic to $\mathbb{S}^1\times [0,1)$ on which $\phi_p$ defines a $C^{\infty}$ diffeomorphism $\phi_p^*$. 

Trivially, the diffeomorphism $\phi_p^*$ restricted on $\mathbb{S}^1\times (0,1)$ is smoothly conjugated to $\phi_p$ restricted on $\mathbb{R}^2\setminus\{0\}$. More specifically, there exists a $C^{\infty}$ map $\Pi: \mathbb{S}^1\times [0,1)\rightarrow \mathbb{R}^2$ such that $\Pi(\mathbb{S}^1\times \{0\})=0$, $\Pi_{|\mathbb{S}^1\times (0,1)}: \mathbb{S}^1\times (0,1)\rightarrow \mathbb{R}^2-\{0\}$ is a diffeomorphism and $$\Pi\circ\phi^*_p = \phi_p\circ \Pi$$ 

\begin{figure}
    \centering
    \includegraphics[scale=0.12]{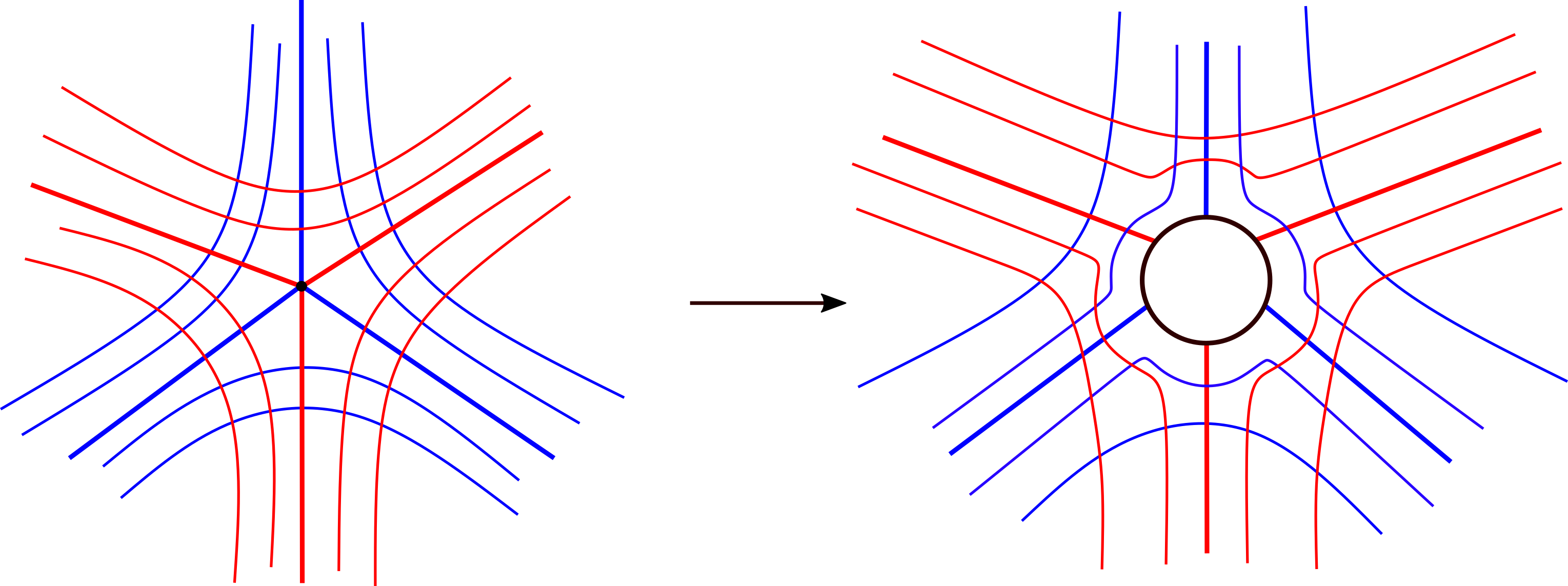}
    \caption{Foliations before and after blow-up}
    \label{f.folibeforeafter}
\end{figure}

Let $(\mathcal{F}^{s,u}_{p})^*:=\{\text{Clos}\big(\Pi^{-1}(L-0_{\mathbb{R}^2})\big)|L\in\mathcal{F}^{s,u}_{p}\}$ be the \emph{stable and unstable foliations of $\phi_p^*$}, where $\text{Clos}(A)$ denotes the closure of a set $A$. It is not hard to prove the following :
\begin{rema}\label{r.phistar}\quad
\newline{\quad}
\begin{itemize}
    \item $(\mathcal{F}^s_{p})^*,(\mathcal{F}^u_{p})^*$ are invariant by $\phi^*_p$ and form a pair of transverse foliations on $\mathbb{S}^1\times (0,1)$ (see Figure \ref{f.folibeforeafter})
    \item $\phi_p^*$ is a Morse-Smale diffeomorphism on $\mathbb{S}^1\times \{0\}$ with $2p$ hyperbolic fixed points ($p$ attracting and $p$ repelling). Each of the previous attracting (resp. repelling) fixed points is contained in a unique leaf of $(\mathcal{F}^u_{p})^*$ (resp. $(\mathcal{F}^s_{p})^*$), whose image by $\Pi$ is an unstable (resp. stable) prong of the origin in $(\mathbb{R}^2, \mathcal{F}^s_{p},\mathcal{F}^u_{p})$ 
\end{itemize}
\end{rema}

Let $R_{\theta}^*$ be the rotation of angle $2\pi\theta$ on $\mathbb{S}^1\times [0,1)$, or equivalently the unique diffeomorphism on $\mathbb{S}^1\times [0,1)$ such that $\Pi\circ R_{\theta}^*=R_{\theta}\circ \Pi$. We define for every $p\geq 1$ and $k\in \llbracket 0, p-1\rrbracket$ $$\phi^*_{pk}:=\phi^*_p\circ R_{k/p}^*= R_{k/p}^*\circ\phi^*_p$$ Consider now the smooth manifold defined as follows: 
$$N^*_{pk}:= \frac{\mathbb{S}^1\times [0,1)\times \mathbb{R}}{((x,y),t+1) \sim (\phi_{pk}^*(x,y),t) }$$
endowed with the constant speed vertical flow $\Phi_{pk}^*=\big((\Phi^*_{pk})^t\big)_{t\in \mathbb{R}}$ given by the vector field $\frac{\partial}{\partial t}$. We will call $(N^*_{pk},\Phi_{pk}^*) $ the \emph{flow obtained from $\Phi_{pk}$ after blowing up $\gamma_{pk}$}. 

Just as in the case of $\phi_p$ and $\phi_p^*$, the flow $\Phi_{pk}^*$ restricted on the interior of $N_{pk}^*$ is $C^{\infty}$ conjugated to $\Phi_{pk}$ restricted on $N_{pk}-\gamma_{pk}$. More specifically, thanks to $\Pi$, we can define a $C^{\infty}$ map $\Pi^*: N^*_{pk}\rightarrow N_{pk}$ such that $\Pi^*_{|\inte{N_{pk}^*}}: \inte{N_{pk}^*}\rightarrow N_{pk}-\gamma_{pk}$ is a $C^{\infty}$ diffeomorphism and for every $t\in \mathbb{R}$ 
\begin{equation} \label{eq.conjugationflowsblowup}
    \Pi^*\circ(\Phi^*_{pk})^t= \Phi_{pk}^t\circ \Pi^*
\end{equation}
We will call $\Pi^*$ the \emph{blow-down map associated to the blow-up of $\gamma_{pk}$}. Moreover, similarly to $\Phi_{pk}$, the flow $\Phi^*_{pk}$ preserves the projections of $(\mathcal{F}^s_p)^*\times \mathbb{R}, (\mathcal{F}^u_p)^*\times \mathbb{R}$ on $N^*_{pk}$ which we will denote respectively by $(F^s_{pk})^*,(F^u_{pk})^*$ and we will call the \emph{stable and unstable foliations of $\Phi_{pk}^*$}. Finally, by our construction and thanks to Remark \ref{r.phistar}, it is easy to see that: 
\begin{rema}\label{r.conjugation1}
\quad \quad
\begin{itemize}
    \item the flow $\Phi^*_{pk}$ restricted on the boundary of $N^*_{pk}$ is a Morse-Smale flow with $gcd(k,p)$ (gcd stands for the greatest common divisor) attracting and $gcd(k,p)$ repelling hyperbolic periodic orbits
    \item $(F^s_{pk})^*,(F^u_{pk})^*$ define a pair of transverse foliations in the interior of $N^*_{pk}$
    \item $(F^{s,u}_{pk})^*:=\{\text{Clos}\big((\Pi^*)^{-1}(L-\gamma_{pk})\big)|L\in F^{s,u}_{pk}\}$
    %\item If $V_{pk}$ is a compact neighborhood of $\gamma_{pk}$ and $d_{pk}$ is any distance on $V_{pk}$ given by a Riemannian metric, then for any $x\in \gamma_{pk}$ and any $y\in F^s_{pk}(x)$ (resp. $y\in F_{pk}^u(x)$) there exists an increasing homeomorphism $h:\mathbb{R}\rightarrow \mathbb{R}$ such that $$d_{pk}((\Phi_{pk})^t(x),(\Phi_{pk})^{h(t)}(y))\underset{t\rightarrow +\infty}{\longrightarrow} 0~~ \big(\text{resp. }d_{pk}((\Phi_{pk})^t(x),(\Phi_{pk})^{h(t)}(y))\underset{t\rightarrow -\infty}{\longrightarrow}0\big) $$
    %\item If $V^*_{pk}$ is a compact neighborhood of $\partial N_{pk}^*$,  $d^*_{pk}$ is any distance on $V^*_{pk}$ given by a Riemannian metric and  $x\in\partial N_{pk}^* $ belongs to a periodic orbit of $\Phi^*_{pk}$, then for any $y\in (F^s_{pk})^*(x)$ (resp. $y\in (F_{pk}^u)^*(x)$) there exists an increasing homeomorphism $h:\mathbb{R}\rightarrow \mathbb{R}$ such that $$d^*_{pk}((\Phi^*_{pk})^t(x),(\Phi^*_{pk})^{h(t)}(y))\underset{t\rightarrow +\infty}{\longrightarrow} 0~~ \big(\text{resp. }d^*_{pk}((\Phi^*_{pk})^t(x),(\Phi^*_{pk})^{h(t)}(y))\underset{t\rightarrow -\infty}{\longrightarrow}0\big) $$
\end{itemize}
\end{rema}
\subsection{Blow-up operation for almost pseudo-Anosov flows}\label{s.blowupflows}

Consider $M$ a closed 3-manifold, $\Phi=(\Phi^t)_{t\in \mathbb{R}}$ an almost pseudo-Anosov flow on $M$, $F^s,F^u$ its stable and unstable foliations and $\gamma$ a periodic orbit of $\Phi$, whose tubular neighborhoods are solid tori. 

Let us mention here that for any pseudo-Anosov (or equivalently expansive) flow $X$ on $M$ such periodic orbits are quite common. Indeed, by Proposition 6 of \cite{Markovnontransitive}, periodic orbits form a dense set inside the non-wandering set of any expansive flow. If $M$ is orientable, the tubular neighborhoods of any periodic orbit of $X$ are solid tori. If $M$ is not orientable, one can prove using Markov partitions that the set of periodic orbits whose tubular neighborhoods are solid tori is dense in the non-wandering set of $X$. 

Consider $(N_{pk},\Phi_{pk})$ the local model for a pseudo-hyperbolic periodic orbit $\gamma_{pk}$ with $p\geq 2$ prongs and  rotation $k \in \llbracket 0, p-1 \rrbracket$. By Proposition \ref{p.aroundcircleprong}, there exist $p\geq 2$, $k \in \llbracket 0, p-1 \rrbracket$, $U_{\gamma}$ a neighborhood of $\gamma$ homeomorphic to a solid torus, $V_{pk}$ a neighborhood of $\gamma_{pk}$ and $H: U_{\gamma}\rightarrow V_{pk}$ a homeomorphism defining an orbital equivalence between $(\Phi,U_\gamma)$ and $(\Phi_{pk},V_{pk})$ and sending stable/unstable leaves in  $U_\gamma$ to stable/usntable leaves in $V_{pk}$. 

Let $(N^*_{pk},\Phi^*_{pk})$ be the flow obtained from $\Phi_{pk}$ by blowing up $\gamma_{pk}$, $\Pi^*:N^*_{pk}\rightarrow N_{pk}$ be the blow-down map associated to the blow-up of $\gamma_{pk}$ and $V^*_{pk}:=(\Pi^*)^{-1}(V_{pk})$. Thanks to Equation \ref{eq.conjugationflowsblowup}, the map $H^{-1}\circ\Pi^*$ defines an orbital equivalence between the flows $(\Phi_{pk}^*,V^*_{pk}-\partial N^*_{pk})$ and $(\Phi,U_\gamma-\gamma)$. Furthermore, for any $x\in V^*_{pk}$ the map $H^{-1}\circ\Pi^*$ takes the connected component of $(F^{s,u}_{pk})^*(x)\cap (V^*_{pk}-\partial N^*_{pk})$ containing $x$ to the connected component of $F^s(H^{-1}\circ\Pi^*(x))\cap (U_{\gamma}-\gamma)$ containing $H^{-1}\circ\Pi^*(x)$.  

Consider $M^*$ the $C^0$ manifold obtained by glueing $V^*_{pk}$ and $M-\inte{U_\gamma}$ along $(\Pi^*)^{-1}(\partial V_{pk})$ and $\partial U_\gamma$ via the map $H^{-1}\circ\Pi^*$. We will call $M^*$ the \emph{manifold obtained after blowing up $\gamma$}. By using the previous properties of $H^{-1}\circ\Pi^*$, we get that the flows $(\Phi_{pk}^*,V^*_{pk})$ and $(\Phi,M-U_\gamma)$ define on $M^*$ a $C^0$ orientable foliation of dimension $1$ that can be parametrized by a flow $\Phi^*=((\Phi^*)^t)_{t\in \mathbb{R}}$, thanks to Theorem 27A of \cite{Whitney}. We will call $\Phi^*$ a \emph{flow obtained after blowing up $\gamma$}. Similarly, the foliations $(F^{s,u}_{pk})^*$ on $ V^*_{pk}$ and $F^{s,u}$ on $M-U_\gamma$ define on $M^*$ two sets of leaves, that will be denoted by $(F^s)^*$ and  $(F^u)^*$ and will be called respectively the \emph{stable and unstable foliations associated to $(M^*,\Phi^*)$}. It is not difficult to see, thanks to Remark \ref{r.conjugation1}, that: 
\begin{rema}\label{r.conjugation2}
\quad \quad 
    \begin{enumerate}
        \item $(F^{s,u})^*$ define a pair of transverse singular foliations in the interior of $M^*$, which are invariant by $\Phi^*$
        \item $gcd(k,p)$ leaves of $(F^s)^*$ and $gcd(k,p)$ leaves of $(F^u)^*$ intersect $\partial M^*$. Each of the previous intersections corresponds to a unique periodic orbit of $\Phi^*$ on $\partial M^*$
        \item contrary to $\Phi^*_{pk}$, even if we fix our choices of $U_\gamma$ and $H$, a flow obtained after blowing up $\gamma$ is uniquely defined up to reparametrization
    \end{enumerate}
\end{rema}

As in the case of hyperbolic models, we can define a blow-down map for $M^*$. Consider $\Pi^*_M: M^*\rightarrow M$ the map defined as the identity on $M-U_\gamma\subset M^*$ and as $H^{-1}\circ \Pi^*$ on $ V^*_{pk}\subset M^*$. This map is trivially continuous and will be called the \emph{blow-down map associated to $(M^*,\Phi^*)$}. By construction, 
\begin{rema}\label{r.conjugationpstar}
\begin{enumerate}
 \item $(F^{s,u})^*=\{\text{Clos}\big((\Pi_M^*)^{-1}(L-\gamma)\big)|L\in F^{s,u}\}$ \item $\Pi^*_M$ defines an orbital equivalence between $(\inte{M^*},\Phi^*)$ and  $(M-\gamma,\Phi)$
 \end{enumerate}
\end{rema}
Finally, as in the case of $F^{s,u}$ we have that:
\begin{prop}\label{p.stablestar}For any metric $d^*$ on $M^*$ (compatible with the the topology of $M^*$), any $x^*\in M^*$ and any $y^*\in (F^s)^*(x)$ (resp. $y^*\in (F^u)^*(x)$), there exists an increasing homeomorphism $h:\mathbb{R}\rightarrow \mathbb{R}$ such that 
\begin{equation*} d^*((\Phi^*)^t(x^*),(\Phi^*)^{h(t)}(y^*))\underset{t\rightarrow +\infty}{\longrightarrow} 0~~ \big(\text{resp. }d^*((\Phi^*)^t(x^*),(\Phi^*)^{h(t)}(y^*))\underset{t\rightarrow -\infty}{\longrightarrow}0\big)
\end{equation*}
\end{prop}
The proof of the above proposition, being a bit long and also being 
intimately related with the proof of our main result, will be postponed to Section \ref{s.proofproposition}.

\subsection{Dehn-Fried surgery for expansive flows}\label{s.dehn-frieddefi}

Following the notations of the previous section, denote by $\Phi^*_{|\partial M^*}$ the restriction of $\Phi^*$ on $\partial M^*$. We will say that a foliation $\mathcal{F}$ on $\partial M^*$ is \emph{good} if it is invariant by $\Phi^*_{|\partial M^*}$ and every leaf of $\mathcal{F}$ is a \emph{global section} of $\Phi^*_{|\partial M^*}$ (i.e. every leaf $l$ of $\mathcal{F}$ is topologically transverse to $\Phi^*_{|\partial M^*}$ and intersects every orbit of $\Phi^*_{|\partial M^*}$) on which the orbits of $\Phi^*_{|\partial M^*}$ return in a uniformly bounded time. If such a foliation exists for $\Phi^*$, by crushing its leaves to points, we are going to define a new blow-down operation giving rise to a new flow on a closed manifold. We will call this operation a Dehn-Fried surgery. 

Let us first prove that, up to reparametrizing $\Phi^*$ close to $\partial M^*$, we can always find good foliations on $\partial M^*$. Fix an orientation on the torus $\partial M^*$ and a point $x_0$ belonging to a periodic orbit $\mathcal{O}$ of $\Phi^*$ inside $\partial M^*$. Denote by $ \sigma_0\in \pi_1(\partial M^*,x_0)=\pi_1(\mathbb{T}^2)$ the homotopy class of $\mathcal{O}$ (endowed with its dynamical orientation). 

\begin{lemm}\label{l.goodfoli}
   For every indivisible element $\sigma\in \pi_1(\partial M^*,x_0)$ such that $\sigma \neq \pm \sigma_0$, there exists $\widetilde{\Phi}^*$, a flow obtained by reparametrizing $\Phi^*$ close to $\partial M^*$, such that $\widetilde{\Phi}^*$ admits a good foliation $\mathcal{F}_{\sigma}$ on $\partial M^*$, whose every leaf is freely homotopic to $\sigma$. 
\end{lemm}
\begin{proof}
Consider $\delta$ a simple closed curve in $\partial M^*$, whose homotopy class is $\sigma$ (such a simple closed representative exists because $\sigma$ is indivisible). Since $\sigma\neq \pm\sigma_0$, the curve $\delta$ intersects every periodic orbit of $\Phi^*_{|\partial M^*}$. Recall that $\Phi^*_{|\partial M^*}$ is orbitally equivalent to a Morse-Smale flow and more specifically to the suspension of a Morse-Smale diffeomorphism of the circle. Using this fact, one can prove that : 
\begin{enumerate}
    \item $\delta$ can be taken transverse to $\Phi^*_{|\partial M^*}$
    \item assuming that (1) is true, $\delta$ intersects every orbit of $\Phi^*_{|\partial M^*}$ in a uniformly bounded time. In this way,  $\delta$ becomes a global section of $\Phi^*_{|\partial M^*}$ with a uniformly bounded return time
    \item the lift of $\Phi^*_{|\partial M^*}$ on $\mathbb{R}^2$, the universal cover of $\partial M^*$, admits a global section intersecting once every orbit of the lift. Thanks to this, the lift of $\Phi^*_{|\partial M^*}$ on $\mathbb{R}^2$ is conjugated to a constant speed vertical flow on $\mathbb{R}^2$. It follows that by reparametrizing $\Phi^*$ close to $\partial M^*$, we can obtain a flow $\widetilde{\Phi}^*=\big((\widetilde{\Phi}^t)^*\big)_{t\in \mathbb{R}}$ for which $(\widetilde{\Phi}^1)^*(\delta)=\delta$
\end{enumerate}
Finally, by pushing $\delta$ by $\widetilde{\Phi}^*$, we obtain a foliation on $\partial M^*$ with the desired properties. 
\end{proof}

Consider $\widetilde{\Phi}^*$ a flow obtained by reparametrizing $\Phi^*$ inside $(\Pi_M^*)^{-1}(U_\gamma)$ and $\mathcal{F}$ a good foliation of $\widetilde{\Phi}^*$, whose leaves, up to free homotopy, belong to the indivisible class $\sigma\neq \pm\sigma_0\in \pi_1(\partial M^*,x_0)$. Let $M_{\mathcal{F}}:=\quotient{M^*}{\mathcal{F}}$ be the manifold obtained from $M^*$ by identifying every leaf of $\mathcal{F}$ to a point. Consider $\Pi_{\mathcal{F}}$ the natural projection from $M^*$ to $M_{\mathcal{F}}$. The map $\Pi_{\mathcal{F}}$ is clearly continuous and will be called the \emph{blow-down map associated to $\mathcal{F}$}. Since $\mathcal{F}$ is $\widetilde{\Phi}^*$ invariant, the flow $\widetilde{\Phi}^*$ projects to a $C^0$ flow $\Phi_{\mathcal{F}}=(\Phi_{\mathcal{F}}^t)_{t\in \mathbb{R}}$ on $M_{\mathcal{F}}$ for which $$\Pi_\mathcal{F}\circ(\widetilde{\Phi}^*)^t= \Phi_{\mathcal{F}}^t\circ \Pi_{\mathcal{F}}$$

Recall that (see Remark \ref{r.conjugation2}) a flow obtained after blowing up $\gamma$ is not uniquely defined and notice that $\Phi_{\mathcal{F}}$ was constructed thanks to a choice of $\widetilde{\Phi}^*$ and $\mathcal{F}$. Despite the previous  facts, it is not difficult to prove that
\begin{prop}
    The orbital equivalence class of $(M_{\mathcal{F}},\Phi_{\mathcal{F}})$ depends only on $\sigma$ and not our choice of $\Phi^*$, $\widetilde{\Phi}^*$ or $\mathcal{F}$.
\end{prop}
This justifies the following definition: 
\begin{defi}\label{d.defisurgery}
    We will say that the flow $(M_{\mathcal{F}},\Phi_{\mathcal{F}})$ is obtained from $(M,\Phi)$ by the \emph{Dehn-Fried surgery along $\gamma$ associated to the homotopy class $\sigma$}.
\end{defi}
\begin{rema}Dehn-Fried surgeries are only defined for periodic orbits with orientable  tubular neighborhoods (or equivalently tubular neighborhoods homeomorphic to a solid torus). The reason for this is that there exists a unique (up to homotopy) foliation on the Klein bottle with no one sided leaf. 
\end{rema}
By construction, $\gamma_{\mathcal{F}}:=\Pi_{\mathcal{F}}(\partial M^*)$ is a periodic orbit of the flow $\Phi_{\mathcal{F}}$ whose tubular neighborhoods are solid tori and also 
\begin{rema}\label{r.conjugationpif}
$\Pi_{\mathcal{F}}$ defines an orbital equivalence between the flows $(\inte{M^*},\Phi^*)$ and $(M_{\mathcal{F}}-\gamma_{\mathcal{F}},\Phi_{\mathcal{F}})$
\end{rema}

Next, just as in the case of $\Phi^*$ or $\widetilde{\Phi}^*$, the flow $\Phi_{\mathcal{F}}$ admits stable and unstable foliations.  Let $(L^{s,u})^*$ be the union of all stable/unstable leaves in $(F^{s,u})^*$ that intersect $\partial M^*$ and 
$$F^{s,u}_{\mathcal{F}}:= \{\Pi_{\mathcal{F}}(L)|L\in(F^{s,u})^*, L\cap \partial M^*=\emptyset\}\cup \{\Pi_{\mathcal{F}}\big((L^{s,u})^*\big)\}$$
\begin{prdef}\label{p.foliationsaftersurgery}$F^{s,u}_{\mathcal{F}}$ form a pair of transverse singular foliations, called respectively the \emph{stable and unstable foliations of $\Phi_{\mathcal{F}}$}. 
\end{prdef}
\begin{proof}
     Since $\Pi_{\mathcal{F}}$ defines a homeomorphism between $\inte{M^*}$ and $M_{\mathcal{F}}-\gamma_{\mathcal{F}}$, by Remark \ref{r.conjugation2} we get that $F^{s,u}_{\mathcal{F}}$ form a pair of transverse singular foliations on $M_{\mathcal{F}}-\gamma_{\mathcal{F}}$. 

    It suffices to show that either $F^{s,u}_{\mathcal{F}}$ form a pair of transverse foliations near $\gamma_{\mathcal{F}}$ or that $\gamma_{\mathcal{F}}$ is a circle prong singularity of $F^{s,u}_{\mathcal{F}}$. Indeed, consider $R$ a transverse standard polygon in $M$ containing a point $x_\gamma\in \gamma$ in its interior. Let $R_1,R_2,...,R_{2p}$ be the closures of the quadrants of $x_\gamma$ in $R$. Notice that $R_1,...,R_n$ are transverse standard rectangles in $M$ with a corner point in $\gamma$. 
    \begin{figure}[h!]
        \centering
        \includegraphics[scale=0.13]{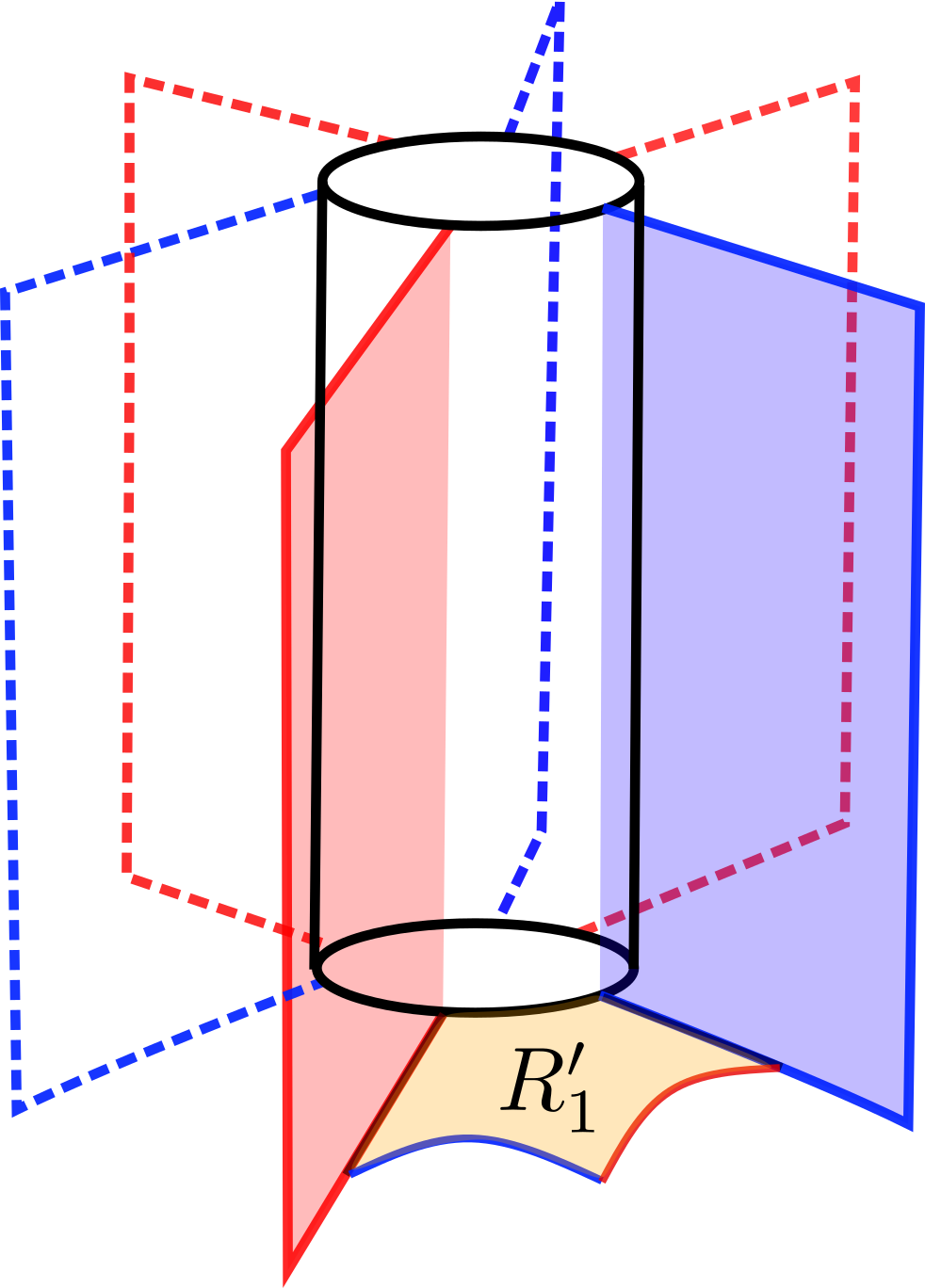}
        \caption{}
        \label{f.liftrectangle}
    \end{figure}

    Denote by $R'_1$ the closure of $(\Pi^*_M)^{-1}(R_1-\gamma)$ in $M^*$. Since $(F^{s,u})^*=\{\text{Clos}\big((\Pi_M^*)^{-1}(L-\gamma)\big)|L\in F^{s,u}\}$, we have that $R'_1$ is homeomorphic to a closed disk that is bounded by two stable segments, two unstable segments and one arc in $\partial M^*$ as in Figure \ref{f.liftrectangle}. We will call each of the previous stable (resp. unstable) segments a \emph{stable} (resp. \emph{unstable}) \emph{boundary component of $R'_1$}. Fix $L$ a leaf of $\mathcal{F}$. By pushing $R'_1$ along the flowlines of $\widetilde{\Phi^*}$, we can construct a topological disk $R^*_1$ transverse to $\widetilde{\Phi}$, bounded by two stable segments, two unstable segments and an arc in  $\partial M^*$ that is included in $L$. By lifting in a similar way $R_2$ on $M^*$, we can define $R^*_2$ such that $R_1^*$ and $R_2^*$ have a common unstable boundary component and $R^*_2\cap \partial M^* \subset L$. By repeating this procedure a finite number of times (eventually by restarting from $R_1$, once we reach $R_{2p}$), we obtain $R_1^*,R_2^*,...,R_{2n}^*$ such that 
    \begin{itemize}
        \item $R_{2k-1}^*$ and $R_{2k}^*$ have a common unstable boundary component and no other intersection for every $k\in\llbracket 1, n\rrbracket$ 
        \item $R_{2k}^*$ and $R_{2k+1}^*$ have a common stable boundary component and no other intersection for every $k\in\llbracket 1, n-1\rrbracket$
        \item $\underset{i\in\llbracket 1, n\rrbracket}{\cup}R_i^*\cap \partial M^*= L$ 
    \end{itemize}
Notice that the number $n$ does not depend on our choice of $L$. Indeed, one can easily prove that $n=K\times gcd(k,p)$, where $K$ is the number of intersections between any leaf of $\mathcal{F}$ and any periodic orbit of $\Phi^*$ on $\partial M^*$ and $gcd(k,p)$ is the number of attracting or repelling periodic orbits of $\widetilde{\Phi^*}$ on $\partial M^*$.

By cutting the disk $R_1^*$ or $R_{2n}^*$ and pushing $R_{2n}^*$ along the flowlines of $\widetilde{\Phi^*}$, in addition to the previous properties, we can also assume that $R_1^*$ and $R_{2n}^*$ have a common stable boundary component and no other intersection. Moreover, by eventually choosing a smaller transverse polygon $R$ and pushing the lifts $R_1^*,R_2^*,...,R_{2n}^*$ along the flowlines of $\widetilde{\Phi^*}$, we can assume without any loss of generality that $R_i^*\cap R_j^*=\emptyset$ for any two distinct integers $i,j$ in $\llbracket 1, 2n\rrbracket$, which are non-consecutive modulo $2n$. 

Thanks to our previous discussion, $\mathcal{A}:=\underset{i\in\llbracket 1, 2n\rrbracket}{\cup}R_i^*$ is an embedded annulus in $M^*$ that is topologically transverse to $\widetilde{\Phi^*}$ and for which $L=\mathcal{A}\cap \partial M^*$. Denote by $x_L\in \gamma_{\mathcal{F}}$ the point $\Pi_{\mathcal{F}}(L)$ and let  $\widetilde{p}:=K\cdot gcd(k,p)$. By construction, $\Pi_{\mathcal{F}}(\mathcal{A})$ is homeomorphic to a closed disk containing $x_L$ and consisting of $2\widetilde{p}$ rectangles in $M_\mathcal{F}$ that are topologically transverse to $\Phi_{\mathcal{F}}$. Thanks to this fact, following the notations of Definition \ref{d.folisingular}, for any sufficiently small $\epsilon>0$, there exists a homeomorphism  $h:\underset{t\in(-\epsilon, \epsilon)}{\cup}\Phi_{\mathcal{F}}^t (\Pi_{\mathcal{F}}(\mathcal{A}))\rightarrow \mathcal{D}_{\widetilde{p}}\times[0,1]$ such that $h(x_L)=(0,\frac{1}{2})$ and $$h\big(\underset{t\in(-\epsilon, \epsilon)}{\cup}\Phi_\mathcal{F}^t (\Pi_{\mathcal{F}}(\mathcal{A}))\cap F^{s,u}_\mathcal{F}\big)=\mathcal{F}_{\widetilde{p}}^{s,u}$$ 
\end{proof}
Thanks to Remarks \ref{r.conjugation1} and \ref{r.conjugationpstar},  $\Pi_{\mathcal{F}} \circ (\Pi^*)^{-1}$ defines a homeomorphism from  $M-\gamma$ to $M_\mathcal{F}-\gamma_F$ that carries stable/unstable leaves of $\Phi$ to stable/unstable leaves of $\Phi_\mathcal{F}$. It follows that the singular foliations $F^{s,u}$ and $F^{s,u}_\mathcal{F}$ are the same (up to homeomorphism) on $M-\gamma$ and $M_\mathcal{F}-\gamma_{\mathcal{F}}$ respectively. However, thanks to our above proof, their behaviors ``at infinity" may differ significantly: 
\begin{rema}\label{r.numberprongs}
    If $K$ is the number of intersections between any leaf of $\mathcal{F}$ and any periodic orbit of $\widetilde{\Phi^*}$ on $\partial{M^*}$, then $\gamma_{\mathcal{F}}$ is a circle $K\times gcd(k,p)$-prong singularity of $F^{s,u}_{\mathcal{F}}$. 
\end{rema} 
Thanks to the continuity of $\Pi_\mathcal{F}$ and Proposition \ref{p.stablestar}, we get that for any metric $d_\mathcal{F}$ on $M_\mathcal{F}$ (compatible with the topology of $M_\mathcal{F}$), any $x\in M_\mathcal{F}$ and any $y\in F^s_\mathcal{F}(x)$ (resp. $y\in F^u_\mathcal{F}(x)$), there exists an increasing homeomorphism $h:\mathbb{R}\rightarrow \mathbb{R}$ such that 
\begin{equation} d_\mathcal{F}(\Phi_\mathcal{F}^t(x),\Phi_\mathcal{F}^{h(t)}(y))\underset{t\rightarrow +\infty}{\longrightarrow} 0~~ \big(\text{resp. }d_\mathcal{F}(\Phi_\mathcal{F}^t(x),\Phi_\mathcal{F}^{h(t)}(y))\underset{t\rightarrow -\infty}{\longrightarrow}0\big)
\end{equation}

Thanks to this fact and the fact that for any point $x\in \gamma_\mathcal{F}$ there exists a transverse standard polygon in $M_\mathcal{F}$ containing $x$ in its interior (this is a result of our proof of Proposition-Definition \ref{p.foliationsaftersurgery} and the fact that $F^{s,u}$ and $F^{s,u}_\mathcal{F}$ coincide on $M-\gamma$ and $M_\mathcal{F}-\gamma_{\mathcal{F}}$), we get that Proposition \ref{p.aroundcircleprong} remains true for $\Phi_\mathcal{F}$. More specifically, 
\begin{prop}\label{p.aroundcircleprongaftersurgery}
    Consider $(N_{p'k'},\Phi_{p'k'})$ the local model for a pseudo-hyperbolic periodic orbit $\gamma_{p'k'}$ with $p'\geq 1$ prongs and  rotation $k' \in \llbracket 0, p'-1 \rrbracket$. Fix $\gamma'$ a periodic orbit of $\Phi_\mathcal{F}$, whose tubular neighborhoods are solid tori. There exist $p'\geq 2$, $k' \in \llbracket 0, p'-1 \rrbracket$, $U_{\gamma'}$ a neighborhood of $\gamma'$ homeomorphic to a solid torus, $V_{p'k'}$ a neighborhood of $\gamma_{p'k'}$ and a homeomorphism $H':U_{\gamma'}\rightarrow V_{p'k'}$ such that 
    \begin{enumerate}
    \item $H'$ defines an orbital equivalence between  $(U_{\gamma'},\Phi_\mathcal{F})$ and  $(V_{p'k'},\Phi_{p'k'})$ 
        \item for any $x\in U_{\gamma'}$ the map $H'$ takes the connected component of $F^s_\mathcal{F}(x)\cap U_{\gamma_\mathcal{F}}$ containing $x$ to the connected component of $F^s_{p'k'}(H'(x))\cap V_{p'k'}$ containing $H'(x)$  
    \end{enumerate}
\end{prop}

To conclude, by our previous arguments, the flow $\Phi_\mathcal{F}$ obtained by a Dehn-Fried surgery on an almost pseudo-Anosov flow $\Phi$ is also almost pseudo-Anosov if and only if  $F^{s,u}_\mathcal{F}$ do not admit an $1$-prong singularity. In simpler words, Dehn-Fried surgeries preserve the class of almost pseudo-Anosov flows, as long as no $1$-prong singularities are produced. Our goal from now on will consist in showing that under the same hypothesis, Dehn-Fried surgeries preserve the class of expansive flows.

\section{Proof of the main result}

\subsection{Cancelling the effect of a Dehn-Fried surgery}\label{ss.cancelsurgery}
Consider $M$ a closed 3-manifold, $\Phi=(\Phi^t)_{t\in \mathbb{R}}$ an almost pseudo-Anosov flow on $M$ and $\gamma$ a periodic orbit of $\Phi$, whose tubular neighborhoods are solid tori. Denote by $(M^*, \Phi^*)$ a flow obtained after blowing up $\gamma$ and by $\Pi^*_M$ its associated blow-down map.

By eventually reparametrizing $\Phi^*$ close to $\partial M^*$, we have that $\Phi^*$ preserves a good foliation $\mathcal{F}$ on $\partial M^*$. Let $M_\mathcal{F}$ be the quotient of $M^*$ by $\mathcal{F}$, $\Pi_\mathcal{F}$ be the projection from $M^*$ to $M_\mathcal{F}$, $\Phi_\mathcal{F}:= \Pi_\mathcal{F}(\Phi^*)$ and  $\gamma_\mathcal{F}:=\Pi_\mathcal{F}(\partial M^*)$. Recall that $\gamma_\mathcal{F}$ is a periodic orbit of $\Phi_\mathcal{F}$. Assume for the rest of this section that  $\Phi_\mathcal{F}$ is an almost pseudo-Anosov flow. Our goal in this small section consists in proving that 

\begin{prop}\label{p.inversesurgery}
    There exists a Dehn-Fried surgery along $\gamma_\mathcal{F}$ that produces a flow orbitally equivalent to $\Phi$.
\end{prop}
\begin{proof}
Let $R$ be a standard transverse polygon in $M$ containing in its interior a point $x_\gamma\in \gamma$. Notice that $R^*:=(\Pi_M^*)^{-1}(R)$ is an embedded annulus in $M^*$ that is transverse to $\Phi^*$ and that $R_\mathcal{F}:=\Pi_\mathcal{F}(R^*)$ is the image of an immersion $\varphi:\mathbb{S}^1\times [0,1] \rightarrow M_\mathcal{F}$ such that $\varphi(\mathbb{S}^1\times \{0\})\subseteq \gamma_\mathcal{F}$ and $\varphi_{|\mathbb{S}^1\times(0,1]}$ is an embedding.

Consider now $(M_\mathcal{F}^*, \Phi_\mathcal{F}^*)$ a flow obtained after blowing up $\gamma_\mathcal{F}$ and $\Pi^*_{M^*}$ its associated blow-down map. 

%Take  $U_{\gamma_\mathcal{F}}$ to be the toric neighborhood of $\gamma_\mathcal{F}$ used in the construction of $M_\mathcal{F}^*$ (see Section \ref{s.blowupflows}). 
%By eventually taking $R$ smaller, we may assume that $R_\mathcal{F}\subset U_{\gamma_{\mathcal{F}}}$.
%Recall that by our construction in Section \ref{s.blowupflows}, there exists $(N_{pk},\Phi_{pk})$ a local model for a pseudo-hyperbolic periodic orbit $\gamma_{pk}$ with $p\geq 2$ prongs and rotation $k \in \llbracket 0, p-1 \rrbracket$, $V_{pk}$ a neighborhood of $\gamma_{pk}$ and $H: U_{\gamma}\rightarrow V_{pk}$ a homeomorphism defining an orbital equivalence between $(\Phi,U_\gamma)$ and $(\Phi_{pk},V_{pk})$. 
%Let $U_{\gamma_\mathcal{F}}$ be the toric neighborhood of $\gamma_\mathcal{F}$ that was removed from $M_\mathcal{F}$ during the construction of $M_\mathcal{F}^*$. 
%We will assume without any loss of generality that $U_{\gamma_\mathcal{F}}\subset\Pi_\mathcal{F}\circ(\Pi^*)^{-1}(\inte{U_\gamma}) $. 
%Let $(M_{pk}^*, \Phi_{pk}^*)$ be the flow obtained after blowing up $\gamma_{pk}$ and $\Pi^*$ the blow-down map associated to the blow-up of $\gamma_{pk}$.  Also, by eventually pushing $R_\mathcal{F}$ along the flowlines of $\Phi_\matcal{F}$, we may assume that the closure of $(\Pi^*)^{-1}(H(R_\mathcal{F}-\gamma_{pk}))$ is an embedded annulus in $M_{pk}^*$ transverse to $ \Phi_{pk}^*$. 
Since $R^*$ is transverse to $\Phi^*$, we have that there exists $\epsilon>0$ such that every orbit segment of $\Phi_\mathcal{F}$ of size $\epsilon$ intersects $R_\mathcal{F}-\gamma_\mathcal{F}$ at most once. Using this result and the fact that $R_\mathcal{F}-\gamma_\mathcal{F}$ is embedded in $M_\mathcal{F}$, by eventually pushing $R_\mathcal{F}$ along the flowlines of $\Phi_\mathcal{F}$, we have that the closure of $(\Pi_\mathcal{F}^*)^{-1}(R_\mathcal{F}-\gamma_\mathcal{F})$ in $M_\mathcal{F}^*$ is an embedded annulus $R_\mathcal{F}^*$ that is topologically transverse to $\Phi_\mathcal{F}^*$. 

Next, using the fact that every orbit of $\Phi^*$ in $\partial M^*$ intersects $R^*\cap \partial M^*$ in a uniformly bounded time, it is not difficult to see that every orbit of $\Phi^*_\mathcal{F}$ on  $\partial M^*_\mathcal{F}$ also intersects $R_\mathcal{F}^*\cap \partial M^*_\mathcal{F}$ in a uniformly bounded time. As in our proof of Lemma \ref{l.goodfoli}, by reparametrizing $\Phi_\mathcal{F}^*$ close to $\partial M^*_\mathcal{F}$, we can construct a flow $\widetilde{\Phi^*_\mathcal{F}}=\big((\widetilde{\Phi^*_\mathcal{F}})^t\big)_{t\in\mathbb{R}}$ for which $(\widetilde{\Phi^*_\mathcal{F}})^1(R_\mathcal{F}^*\cap \partial M^*_\mathcal{F})=R_\mathcal{F}^*\cap \partial M^*_\mathcal{F}$. The previous flow admits a good foliation $\mathcal{G}$ on $\partial M^*_\mathcal{F}$, whose every leaf is freely homotopic to the curve formed by $R_\mathcal{F}^*\cap \partial M^*_\mathcal{F}$. 

Consider $(M_\mathcal{G}, \Phi_\mathcal{G})$, the flow obtained by quotienting $M_\mathcal{F}^*$ by $\mathcal{G}$, $\Pi_\mathcal{G}$ be the projection from $M_\mathcal{F}
^*$ to $M_\mathcal{G}$, $\gamma_\mathcal{G}:=\Pi_\mathcal{G}(\partial M_\mathcal{F}^*)$ and $R_\mathcal{G}:=\Pi_\mathcal{G}(R_\mathcal{F}^*)$. Notice that $R_\mathcal{G}$ is a transverse standard polygon in $M_\mathcal{G}$ containing a unique point of $\gamma_\mathcal{G}$, say $x_{\gamma_\mathcal{G}}$. 
\vspace{0.2cm}

\textbf{Claim: }The flows $ \Phi_\mathcal{G}$ and $\Phi$ are orbitally equivalent.
\vspace{0.2cm}   
    
    By Remarks \ref{r.conjugationpstar} and  \ref{r.conjugationpif}, we have that $H:=\Pi_\mathcal{G}\circ (\Pi_\mathcal{F}^*)^{-1}\circ \Pi_\mathcal{F}\circ(\Pi_M^*)^{-1}$ defines an orbital equivalence between $(M-\gamma, \Phi)$ and $(M_\mathcal{G}-\gamma_\mathcal{G},\Phi_\mathcal{G})$. Furthermore, by Remark \ref{r.numberprongs}, we have that $\gamma$ and $\gamma_\mathcal{G}$ have the same number of stable/unstable prongs. Even more, by our construction of $R_\mathcal{G}$, by eventually pushing $R_\mathcal{G}$ along the flowlines of $\Phi_\mathcal{G}$, we have that 
    \begin{equation}\label{eq.orbitequivalence}
        H(R-\gamma)=R_\mathcal{G}-\gamma_\mathcal{G}
    \end{equation}
    
    Finally, consider $f_R$ and $f_{R_\mathcal{G}}$ the first return maps of $\Phi$ on $R$ and of $\Phi_\mathcal{G}$ on $R_\mathcal{G}$ respectively. Both $f_R$ and $f_{R_\mathcal{G}}$ are well defined and continuous near $x_\gamma$ and $x_{\gamma_\mathcal{G}}$ respectively. Thanks to \ref{eq.orbitequivalence}, $f_R$ and $f_{R_\mathcal{G}}$ are conjugated near $x_\gamma$ and $x_{\gamma_\mathcal{G}}$ respectively. It follows that, by eventually modifying $H$ near $\gamma$, we can extend $H$ to an orbital equivalence between $ \Phi_\mathcal{G}$ and $\Phi$, which finishes the proof of the proposition.  
\end{proof}
\subsection{On the closeness of points before and after the blow-up operation}
The main difficulty in the proof of Theorem \ref{t.maintheorem} is contained in the fact that two orbits that are close before Dehn-Fried surgery are not necessarily close after the surgery. The goal of this section, consists in proving that any two orbits of a hyperbolic model flow that are close for a long time before a blow-up operation remain also close after the blow-up operation (see Lemma \ref{l.comparisonmetricsflows} for a more detailed statement). Define $ \phi_p, \phi_{pk},\mathcal{F}^s_p, \mathcal{F}^u_p$, $N_{pk}, \Phi_{pk}$, $\gamma_{pk}, F^s_{pk}, F^u_{pk}$ for every $p\geq 1$ and $k\in \llbracket 0, p-1\rrbracket $ as in Section \ref{s.prelim}. 

Consider $\mathbb{R}^2\setminus\{0\}$ endowed with the polar coordinates $\displaystyle{\{(r,\theta)_2|r\in \mathbb{R}^+, \theta\in \frac{\mathbb{R}}{2\pi\mathbb{Z}}\}}$. For every $p\in \mathbb{N}^*$, by mutliplying our previous circle parametrization by $2/p$, we obtain a new set of coordinates $\displaystyle{\{(r,\theta)_p|r\in \mathbb{R}^+, \theta\in \frac{\mathbb{R}}{p\pi\mathbb{Z}}\}}$ in $\mathbb{R}^2\setminus\{0\}$. Denote by $d_p^{pol}$ the distance associated to the Riemannian metric $dr^2+d\theta^2$ (in the $(\cdot,\cdot)_p$ coordinates). Denote also by $d^{eucl}_p$ the distance associated to the Riemannian metric $dr^2+r^2d\theta^2=d(r\cos{\theta})^2+d(r\sin{\theta})^2$ (in the  $(\cdot,\cdot)_p$ coordinates). Notice that $d^{eucl}_2$ defines the regular euclidean distance on $\mathbb{R}^2\setminus\{0\}$. It is now not difficult to prove the following properties for $d_p^{eucl}$ and $d_p^{pol}$:

%and let $d_2^{pol}$ be the distance associated to the Riemannian metric $dr^2+d\theta^2$. By quotienting $\mathbb{R}^2\setminus\{0\}$ by the rotation $(r, \theta)_2 \rightarrow (r, \theta + \pi)_2$, we obtain $\mathbb{R}^2\setminus\{0\}$ endowed with a new set of coordinates $\displaystyle{\{(r,\theta)_1|r\in \mathbb{R}^+, \theta\in \frac{\mathbb{R}}{\pi\mathbb{Z}}\}}$. More specifically, for every $r>0$ and every $\displaystyle{\theta\in \frac{\mathbb{R}}{\pi\mathbb{Z}}}$ we have that $\pi_2^{-1}((r,\theta)_1)=\{(r,\theta)_2,(r,\theta+\pi)_2\}$. Moreover, using the fact that $dr^2+d\theta^2$ (in the $(\cdot,\cdot)_2$ coordinates) is invariant by the action of $(r, \theta)_2 \rightarrow (r, \theta + \pi)_2$, by passing to the quotient, we obtain a new Riemannian metric on $\mathbb{R}^2\setminus\{0\}$, namely $dr^2+d\theta^2$ (in the $(\cdot,\cdot)_1$ coordinates). Denote by $d_1^{pol}$ the distance associated to the previous metric. 

%Similarly, for any $p\in \mathbb{N}^*$ we define a set of coordinates $\displaystyle{\{(r,\theta)_p|r\in \mathbb{R}^+, \theta\in \frac{\mathbb{R}}{p\pi\mathbb{Z}}\}}$ in $\mathbb{R}^2\setminus\{0\}$ such that for every $r>0$ and every $\theta\in \frac{\mathbb{R}}{\pi\mathbb{Z}}$ we have $\pi_p^{-1}((r,\theta)_1)=\{(r,\theta+k\pi)_p|k\in\llbracket0,p-1\rrbracket\}$. Moreover, by lifting by $\pi_p$ the metric $dr^2+d\theta^2$ (in the $(\cdot,\cdot)_1$ coordinates) we obtain a new Riemannian metric on $\mathbb{R}^2\setminus\{0\}$ whose associated distance will be denoted by $d_p^{pol}$. 

\begin{prop}\label{p.metricscomparison}
\begin{enumerate}
\item For every $r_p,r_p'\in \mathbb{R}^+$ and $\theta_p,\theta'_p\in[0,p\pi)$ we have that  $$(d_p^{pol})^2((r_p,\theta_p), (r_p',\theta_p'))=|r_p-r_p'|^2+\min_{k\in\{-1,0,1\}}|\theta_p-\theta_p'+kp\pi|^2$$
\item For every $p\geq 1$, there exist $\epsilon>0$ and two sequences of points $(z_n)_{n\in\mathbb{N}},(z'_n)_{n\in\mathbb{N}}$ in $\mathbb{R}\setminus \{0\}$ such that $z_n,z'_n\underset{n\rightarrow +\infty}{\longrightarrow}0$, thus $d^{eucl}_p(z_n, z'_n) \underset{n\rightarrow +\infty}{\longrightarrow}0$, and $d^{pol}_p(z_n, z'_n)>\epsilon$ for every $n\in \mathbb{N}$
\item Every rotation of the form $(r, \theta)_p \rightarrow (r, \theta + k)_p$ preserves both $d^{eucl}_p$ and $d^{pol}_p$
\item For every $p\geq 1$, the metric spaces $(\mathbb{R}^2\setminus\{0\},d^{eucl}_p)$ and $(\mathbb{R}^2\setminus\{0\},d^{pol}_p)$ are not complete. The metric $d^{eucl}_p$ can be extended into a complete metric on $\mathbb{R}^2$, whereas $d^{pol}_p$ can be extended into a complete metric on $\mathbb{S}^1\times[0,1)$, the surface obtained after blowing-up the origin of $\mathbb{R}^2$
\item Let $p,q\geq 1$, $Q_p$ (resp. $Q_q$) be a quadrant of the origin inside $(\mathbb{R}^2, \mathcal{F}^{s}_p, \mathcal{F}^u_p)$ (resp. $(\mathbb{R}^2, \mathcal{F}^{s}_q, \mathcal{F}^u_q)$). There exist $\theta_0,\theta_1\in \mathbb{R}$ such that  $Q_p=\{(r, \theta)_p|\theta\in[\theta_0,\theta_0+\pi/2] \}$ and $Q_q=\{(r, \theta)_q|\theta\in[\theta_1,\theta_1+\pi/2 ]\}$. The map $\varphi:Q_p\rightarrow Q_q$ defined by $\varphi((r,\theta)_p)=(r,\theta+\theta_1-\theta_0)_q$ defines an isometry between $Q_p$ and $Q_q$ endowed with the distances $d^{pol}_p$ and $d^{pol}_q$ respectively. We obtain the same result by replacing $d^{pol}$ by $d^{eucl}$. 

\end{enumerate}
\end{prop}

Even though two points in $\mathbb{R}^2\setminus\{0\}$ that are $d^{eucl}_p$-close are not necessarily $d^{pol}_p$-close (see Item (2) of the above proposition), in the following lines we will prove that two $\phi_p$-orbits that stay for some time $d^{eucl}_p$-close, also stay $d^{pol}_p$-close (see Proposition \ref{p.metricscomp1}). This result constitutes one of the key arguments in the proof of Lemma \ref{l.comparisonmetricsflows}.

Fix $p\geq 1$. Consider $Q$ a quadrant of the origin inside  $(\mathbb{R}^2, \mathcal{F}^{s}_p, \mathcal{F}^u_p)$ and $L^s,L^u$ the stable and unstable prongs of the origin that bound $Q$. For every $x\in Q$ there exists a unique pair $(x_s,x_u)\in L^s\times L^u$ such that $\{x\}=\mathcal{F}_p^s(x_u)\cap \mathcal{F}_p^u(x_s)$. The maps $x\rightarrow x_s$ and $x\rightarrow x_u$ are continuous and will be denoted by $\pi_s$ and $\pi_u$ respectively. Fix $I^s$ and $I^u$ two closed segments in $L^s\setminus\{0_{\mathbb{R}^2}\}$ and $L^u\setminus\{0_{\mathbb{R}^2}\}$ respectively. For every $M,N\in\mathbb{N}$, we will denote by $\mathcal{O}^p_{M,N}(I^s,I^u)$ the set 
of $x\in Q$ such that $\pi_s(\phi_p^{-M}(x))\in I^s$ and $\pi_u(\phi_p^N(x))\in I^u$. Let also $\mathcal{O}^p(I^s,I^u):= \underset{M,N\in\mathbb{N}}{\cup}\mathcal{O}^p_{M,N}(I^s,I^u)$. 
\begin{prop}\label{p.metricscomp1}
    For any $\epsilon>0$, there exists $\eta(\epsilon), \eta'(\epsilon)>0$ with $\eta(\epsilon), \eta'(\epsilon)\underset{\epsilon\rightarrow 0}{\longrightarrow} 0$ such that for any $M,N\in\mathbb{N}$ and every $x,y\in \mathcal{O}_{M,N}(I^s,I^u)$ we have that 
     $$ d^{pol}_p(x,y)< \epsilon\implies  d^{eucl}_p(x,y)< \eta(\epsilon)$$
    $$\max_{k\in\llbracket -M, N\rrbracket}\big(d^{eucl}_p(\phi^k_p(x),\phi^k_p(y))\big)< \epsilon\implies  \max_{k\in\llbracket -M, N\rrbracket}\big(d^{pol}_p(\phi^k_p(x),\phi^k_p(y))\big)< \eta'(\epsilon)$$
   
\end{prop}
\begin{proof}
    Thanks to Item (5) of Proposition \ref{p.metricscomparison}, it suffices to prove the above statement for $p=2$. Notice that $\mathcal{O}^2(I^s,I^u)$ is bounded with respect to both $d^{eucl}_2$ and $d^{pol}_2$. Also, we remark that if $x\in  \mathcal{O}_{M,N}(I^s,I^u)$, then for every $k\in\llbracket -M, N\rrbracket$ we have that $\phi^k_2(x)\in\mathcal{O}^2(I^s,I^u)$. Take $(r_1,\theta_1)_2,(r_2,\theta_2)_2\in \mathcal{O}^2(I^s,I^u)$ such that $d^{pol}_2\big((r_1,\theta_1)_2,(r_2,\theta_2)_2\big)<\epsilon$. 
\begin{equation}\label{eq.euclpolar}
    \begin{aligned}
        (d^{eucl}_2)^2\big((r_1,\theta_1)_2,(r_2,\theta_2)_2\big)&=(r_1\cos{\theta_1}-r_2\cos{\theta_2})^2 + (r_1\sin{\theta_1}-r_2\sin{\theta_2})^2 \\ &=r_1^2+r_2^2-2r_1r_2 \sin{\theta_1}\sin{\theta_2}-2r_1r_2 \cos{\theta_1}\cos{\theta_2} \\&=(r_1-r_2)^2 + 2r_1r_2(1-\cos(\theta_1-\theta_2)) 
    \end{aligned}
 \end{equation}  
    Using Item (1) of Proposition \ref{p.metricscomparison}, the fact that $\mathcal{O}^2(I^s,I^s)$ is bounded and the fact that $d^{pol}_2\big((r_1,\theta_1)_2,(r_2,\theta_2)_2\big)<\epsilon$, we get that there exists $\eta>0$ (depending only on $\epsilon$) such that $ d^{eucl}_2\big((r_1,\theta_1)_2,(r_2,\theta_2)_2\big)<\eta$. This proves the first implication of our proposition.
    
    Let us now prove the second implication. Consider $N,M\in \mathbb{N}$ and $(r_1,\theta_1)_2,(r_2,\theta_2)_2\in \mathcal{O}^2_{M,N}(I^s,I^u)$ such that $\max_{k\in\llbracket -M, N\rrbracket}\big(d^{eucl}_2(\phi^k_2(x),\phi^k_2(y))\big)< \epsilon$. Otherwise said, 
    \begin{equation}\label{eq.hypothesisequiv}
        \max_{k\in\llbracket -M, N\rrbracket}\big[2^{-2k}(r_1\cos{\theta_1}-r_2\cos{\theta_2})^2 + 2^{2k}(r_1\sin{\theta_1}-r_2\sin{\theta_2})^2\big]< \epsilon^2
    \end{equation}

    Suppose without any loss of generality that $N=\min(M,N)$. In this case, by applying \ref{eq.hypothesisequiv} for $k=N$ and $k=-N$, we have that 
$2^{2N}(r_1\cos{\theta_1}-r_2\cos{\theta_2})^2+2^{2N}(r_1\sin{\theta_1}-r_2\sin{\theta_2})^2<2\epsilon^2$. By the same calculations as in \ref{eq.euclpolar}, this implies that 
    \begin{equation}\label{eq.final}
        2^{2N}(r_1-r_2)^2 + 2\cdot 2^{2N}r_1r_2(1-\cos(\theta_1-\theta_2))<2\epsilon^2
    \end{equation}
   
    Consider $a,b,\theta_0\in \mathbb{R}$ such that $I^u=\{(r,\theta)_2| r\in[a,b] \text{ and } \theta=\theta_0 \}$. By definition $\pi_u(\phi_2^N((r_1,\theta_1)_2)),\pi_u(\phi_2^N(r_2,\theta_2)_2)\in I^u$. It follows that $2^Nr_1\geq a$ and $2^Nr_2\geq a$.  Using \ref{eq.final} and Item (1) of Proposition \ref{p.metricscomparison}, we get that there exists $\eta'>0$ (depending only on $\epsilon$) such that $ d^{pol}_2\big((r_1,\theta_1)_2,(r_2,\theta_2)_2\big)<\eta'$. We deduce the second implication of the proposition. 
    %Since both $(r_1,\theta_1)_2$ and $(r_2,\theta_2)_2$ belong in $\mathcal{O}^2_{M,N}(I_h,I_v)$, we have that $\frac{1}{2}|r_2\cos(\theta_2)|\leq|r_1\cos(\theta_1)|\leq 2|r_2\cos(\theta_2)|$ and  $\frac{1}{2}|r_2\sin(\theta_2)|\leq|r_1\sin(\theta_1)|\leq 2|r_2\sin(\theta_2)|$. Therefore, $$\frac{1}{2}\leq\frac{r_1}{r_2}\leq2$$ 
    
   % Next, using Equation \ref{eq.euclpolar}, we get that $|r_1-r_2|<\epsilon$ and by Equation \ref{eq.hypothesisequiv}, we also have that  $2^{N}|r_1\sin{\theta_1}-r_2\sin{\theta_2}|<\epsilon$. Moreover,
    %\begin{align*}
    %      2^{N}|r_1\sin{\theta_1}-r_2\sin{\theta_2}|&\geq 2^N r_1|\sin{\theta_1}-\sin{\theta_2}|-2^N|r_1-r_2||\sin{\theta_1}|//&\geq  
    %\end{align*}
   
    %By Equation \ref{eq.euclpolar}, since the euclidean distance of $(r^n_1,\theta^n_1)_2,(r^n_2,\theta^n_2)_2$ is bounded by $\epsilon_n$, we also have that $|r_1^n-r_2^n|<\epsilon_n$. Therefore, $N_n,M_n\underset{n\rightarrow+\infty}{\longrightarrow} +\infty$. 

    %Assume now that there exist $\eta>0, \epsilon_n\underset{n\rightarrow+\infty}{\longrightarrow}0$, $(r^n_1,\theta^n_1)_2,(r^n_2,\theta^n_2)_2\in Q\setminus\{0\}$, $M_n,N_n\in \mathbb{N}$ such that:
    %\begin{itemize}
       % \item $\pi_h(\phi_2^{-M_n}((r^n_1,\theta^n_1)_2)), \pi_h(\phi_2^{-M_n}((r^n_2,\theta^n_2)_2))\in I_h$
        %\item $\pi_v(\phi_2^{N_n}((r^n_1,\theta^n_1)_2)), \pi_v(\phi_2^{N_n}((r^n_2,\theta^n_2)_2))\in I_v$
        %\item $\max_{k\in\llbracket -M_n, N_n\rrbracket}\big(d^{eucl}_2(\phi^k_2((r^n_1,\theta^n_1)_2),\phi^k_2((r^n_2,\theta^n_2)_2))\big)< \epsilon_n$
        %\item $ d^{pol}_2\big((r^n_1,\theta^n_1)_2,(r^n_2,\theta^n_2)_2\big)>\eta$
   % \end{itemize}

\end{proof}
For any $k\in\llbracket 0, p-1\rrbracket$ consider the rotation $R_{k/p}((r,\theta)_p)=(r,\theta+k\pi)_p$ (or equivalently $R_{k/p}((r,\theta)_2)=(r,\theta+2k\pi/p)_2$). Recall that by definition $\phi_{pk}=\phi_p\circ R_{k/p}=R_{k/p} \circ \phi_p$ and that by Item (3) of Proposition \ref{p.metricscomparison} $R_{k/p}$ is an isometry for both $d^{eucl}_p$ and $d^{pol}_p$. Using the two previous facts, one can easily generalize the previous proposition to the following one. 

Take $k\in\llbracket 0, p-1\rrbracket$. Take $I^s,I^u$ two closed segments in $\mathcal{F}^s_p(0_\mathbb{R}^2)\setminus\{0_\mathbb{R}^2\}$ and $\mathcal{F}^u_p(0_\mathbb{R}^2)\setminus\{0_\mathbb{R}^2\}$ respectively. Define as before the sets $\mathcal{O}_{M,N}(I^s,I^u)$ for every $M,N\in \mathbb{N}$. 
\begin{prop}\label{p.metricscomp2}
     For any $\epsilon>0$, there exists $\eta(\epsilon), \eta'(\epsilon)>0$ with $\eta(\epsilon),\eta'(\epsilon)\underset{\epsilon\rightarrow 0}{\longrightarrow} 0$ such that for any $M,N\in\mathbb{N}$ and every $x,y\in \mathcal{O}_{M,N}(I^s,I^u)$ we have that 
     $$ d^{pol}_p(x,y)< \epsilon\implies  d^{eucl}_p(x,y)< \eta(\epsilon)$$
    $$\max_{l\in\llbracket -M, N\rrbracket}\big(d^{eucl}_p(\phi^l_{pk}(x),\phi^l_{pk}(y))\big)< \epsilon\implies  \max_{l\in\llbracket -M, N\rrbracket}\big(d^{pol}_{pk}(\phi^l_p(x),\phi^l_{pk}(y))\big)< \eta'(\epsilon)$$
   
\end{prop}
Consider from now on $p\geq 2$ and $k\in \llbracket 0, p-1\rrbracket$. Recall that $$N_{pk}:= \frac{\mathbb{R}^2\times \mathbb{R}}{((x,y),t+1) \sim (\phi_{pk}(x,y),t) }$$ 
For any $t\in\mathbb{R}$, denote by $\mathcal{P}_t$ the projection of $\mathbb{R}^2\times \{t\}$ on $N_{pk}$. Notice that $\mathcal{P}_t=\mathcal{P}_{t'}$ if and only if $t-t'\in \mathbb{Z}$. Take $\pi_{t}:\underset{t\in(-1/4,1/4)}{\cup}\Phi_{pk}^t(\mathcal{P}_t)\rightarrow \mathcal{P}_t$ the projection on $\mathcal{P}_t$ along the flowlines of $\Phi_{pk}$. 

Consider $\gamma_{pk}$ the unique periodic orbit of $\Phi_{pk}$, $(N^*_{pk},\Phi^*_{pk})$ the blow-up of $(N_{pk},\Phi_{pk})$ along $\gamma_{pk}$ and $\Pi^*$ its associated blow-down map. Fix $d$ a metric on $N_{pk}$, $d^*$ a metric on $N^*_{pk}$ (compatible with the topologies of $N_{pk}$ and $N^*_{pk}$ respectively) and $V$ a neighborhood of $\gamma_{pk}$ homeomorphic to a solid torus such that for every $t\in [0,1]$ the set $V\cap \mathcal{P}_t$ is a transverse standard polygon for $\Phi_{pk}$. For any $x\in V-\gamma_{pk}$ consider $$T^+(x)=\sup\{t>0|\Phi_{pk}^{t}(x)\in V\}$$ $$T^-(x)=\sup\{t<0|\Phi_{pk}^{t}(x)\in V\}$$ Notice that since $x\notin \gamma_{pk}$ either $T^+(x)$ or $T^-(x)$ is finite. Take $J_{exit}(x)$ to be the closure in $\mathbb{R}$ of the interval $(T^-(x),T^+(x))$.
\begin{lemm}\label{l.comparisonmetricsflows}
    For any $\epsilon>0$ sufficiently small, there exist $\eta''(\epsilon,d,d^*,V)>0$ with $\eta'' \underset{\epsilon\rightarrow 0}{\longrightarrow} 0$ such that for any $x\in V-\gamma_{pk}$, any continuous and strictly increasing map $h:J_{exit}(x)\rightarrow \mathbb{R}$ and any $y\in N_{pk}-\gamma_{pk}$  $$\sup_{t\in J_{exit}(x)}\big(d(\Phi_{pk}^t(x),\Phi_{pk}^{h(t)}(y))\big)< \epsilon\implies  \sup_{t\in J_{exit}(x)}\big(d^*((\Phi_{pk}^*)^t(x^*),((\Phi_{pk}^*)^{h(t)}(y^*)\big)< \eta''$$
    where $x^*:=(\Pi^{*})^{-1}(x)$ and $y^*:=(\Pi^{*})^{-1}(y)$.
    
\end{lemm}
\begin{proof}
    It suffices to prove that above statement for $h=id$ and for $x,y$ inside the same $\mathcal{P}_t$. In other words,

 \textbf{Claim $\text{N}^{o}$ 0:} It suffices to show that for any $\epsilon>0$ sufficiently small there exists $\widetilde{\eta}(\epsilon,d,d^*,V)>0$ with $\widetilde{\eta}\underset{\epsilon\rightarrow 0}{\longrightarrow} 0$, such that for any $s\in \mathbb{R}$, $x\in (V-\gamma_{pk})\cap \mathcal{P}_s$ and any $y\in \mathcal{P}_s-\gamma_{pk}$ $$\sup_{t\in J_{exit}(x)}\big(d(\Phi_{pk}^t(x),\Phi_{pk}^{t}(y))\big)< \epsilon\implies  \sup_{t\in J_{exit}(x)}\big(d^*((\Phi_{pk}^*)^t(x^*),(\Phi_{pk}^*)^{t}(y^*))\big)< \widetilde{\eta}$$
\begin{proof}
    Let us prove Lemma \ref{l.comparisonmetricsflows}, assuming that the above is true. 
    
    First, using the continuity of $\Phi_{pk}$ and the compactness of $V$, one gets that for any $\epsilon>0$ sufficiently small, there exists $\epsilon'(\epsilon)>0$ with $\epsilon' \underset{\epsilon\rightarrow 0}{\longrightarrow} 0$ such that for any $t\in [0,1]$, any $z\in \mathcal{P}_t\cap V$ and any $z'\in N_{pk}$ 
    \begin{equation}\label{eq.uniformprojectionclaim0}
    d(z,z')<\epsilon \implies \pi_t(z')\text{ is well defined and } d(z,\pi_t(z'))<\epsilon'
    \end{equation}

    \vspace{0.2cm}
    Take $x\in V-\gamma_{pk}$, $h:J_{exit}(x)\rightarrow \mathbb{R}$ a continuous and  strictly increasing map and $y\in N_{pk}-\gamma_{pk}$ such that for  every $t\in J_{exit}(x)$ we have $$d(\Phi_{pk}^t(x),\Phi_{pk}^{h(t)}(y))< \epsilon$$ We will assume without any loss of generality that $x\in \mathcal{P}_0$; hence, by definition, $\Phi_{pk}^t(x)\in \mathcal{P}_t$ for every $t\in \mathbb{R}$. Thanks to \ref{eq.uniformprojectionclaim0} and by definition of $\pi_t$, if  $\epsilon$ is sufficiently small, then for every $t\in J_{exit}(x)$ we have that 
    \begin{itemize}
        \item $\pi_{t}(\Phi_{pk}^{h(t)}(y))$ is well defined. Let $Y:=\pi_{0}(\Phi_{pk}^{h(0)}(y))$
        \item $\pi_{t}(\Phi_{pk}^{h(t)}(y))=\Phi_{pk}^t(Y)$
        \item $d(\Phi_{pk}^t(x),\Phi_{pk}^t(Y))< \epsilon'$
    \end{itemize}

   Consider now $t_0\in \mathbb{R}$ close to $h(0)$ such that $Y=\Phi^{t_0}_{pk}(y)$ and $\epsilon''(\epsilon)>0$ such that for every $t\in J_{exit}(x)$ $$d\big(\Phi_{pk}^{h(t)}(y), \pi_{t}(\Phi_{pk}^{h(t)}(y))\big)= d\big(\Phi_{pk}^{h(t)}(y), \Phi_{pk}^{t}(Y)\big)= d\big(\Phi_{pk}^{h(t)}(y), \Phi_{pk}^{t+t_0}(y)\big)<\epsilon''$$ Notice that $\epsilon'' \underset{\epsilon\rightarrow 0}{\longrightarrow} 0$ and that the points  $\Phi_{pk}^{h(t)}(y), \Phi_{pk}^{t+t_0}(y)$ are contained in the $(\epsilon+\epsilon'')$-neighborhood of $V$ for every $t\in J_{exit}(x)$. Thanks to the previous facts, the continuity of $\Phi_{pk}$ and the compactness of the $(\epsilon+\epsilon'')$-neighborhood of $V$, if $\epsilon>0$ (and thus also $\epsilon''$) is sufficiently small, we can find $T(\epsilon)>0$ with $T \underset{\epsilon\rightarrow 0}{\longrightarrow} 0$ such that for every $t\in J_{exit}(x)$ we have $|h(t)-t_0-t|<T$. 
   
   By our initial hypothesis, if $\epsilon$ (and thus also $\epsilon'$) is sufficiently small, this implies that there exists $\widetilde{\eta}(\epsilon,d,d^*,V)>0$ with $\widetilde{\eta} \underset{\epsilon\rightarrow 0}{\longrightarrow} 0$ such that for every $t\in J_{exit}(x)$ we have that 
    \begin{equation}\label{eq.boundafterreparametrization}
    d^*((\Phi_{pk}^*)^t(x^*),(\Phi_{pk}^*)^{t}(Y^*))=d^*((\Phi_{pk}^*)^t(x^*),(\Phi_{pk}^*)^{t+t_0}(y^*))< \widetilde{\eta}
    \end{equation}
    where $x^*:=(\Pi^{*})^{-1}(x)$, $y^*:=(\Pi^{*})^{-1}(y)$ and $Y^*=(\Pi^*)^{-1}(Y)$.
   
   Finally, using the fact that $|h(t)-t_0-t|<T$, the continuity of $\Phi_{pk}^*$, the compactness of the pre-image by $\Pi^*$ of the $(\epsilon+\epsilon'')$-neighborhood of $V$ and Inequality \ref{eq.boundafterreparametrization}, we get that there exists $\eta''(\epsilon,d,d^*,V)>0$ with $\eta'' \underset{\epsilon\rightarrow 0}{\longrightarrow} 0$ such that $$d((\Phi^*_{pk})^t(x^*),(\Phi^*_{pk})^{h(t)}(y^*))< \eta''$$ which proves Lemma \ref{l.comparisonmetricsflows} and thus finishes the proof of the claim. 
\end{proof}
In the previous claim, we showed that it suffices to prove Lemma \ref{l.comparisonmetricsflows} only in the case where $h=id$. We will now show that, when $h=id$, Lemma \ref{l.comparisonmetricsflows} boils down to a discrete time problem. For every $s\in[0,1]$, let $0_{\mathcal{P}_s}$ be the unique point inside $\mathcal{P}_s\cap \gamma_{pk}$.

\textbf{Claim $\text{N}^{o}$ 1:} For $\epsilon$ sufficiently small, it suffices to prove that there exists $\zeta(\epsilon,d,d^*,V)>0$ with $\zeta\underset{\epsilon\rightarrow 0}{\longrightarrow} 0$ such that for every $s\in [0,1]$, every $x\in V\cap (\mathcal{P}_s\setminus 0_{\mathcal{P}_s})$ and every $y\in \mathcal{P}_s\setminus 0_{\mathcal{P}_s}$ $$\sup_{t\in J_{exit}(x)}\big(d(\Phi_{pk}^t(x),\Phi_{pk}^{t}(y))\big)< \epsilon\implies \sup_{n\in \mathbb{Z}\cap J_{exit}(x)}\big(d^{pol}_p(\Phi_{pk}^n(x),\Phi_{pk}^{n}(y))\big)< \zeta$$

 \begin{proof}
      Take $\epsilon$ sufficiently small, $s\in [0,1]$, $x\in  V\cap (\mathcal{P}_s\setminus 0_{\mathcal{P}_s})$ and $y\in \mathcal{P}_s\setminus 0_{\mathcal{P}_s}$ such that $\sup_{t\in J_{exit}(x)}\big(d(\Phi_{pk}^t(x),\Phi_{pk}^{t}(y))\big)<\epsilon$. Consider the metric $d^{pol}_p$ on $\mathcal{P}_s\setminus 0_{\mathcal{P}_s}$. Thanks to Proposition \ref{p.metricscomparison}, by pulling back $d^{pol}_p$ by $\Pi^*$, we can define a metric on $(\Pi^*)^{-1}(\mathcal{P}_s\setminus 0_{\mathcal{P}_s})$, which extends to a complete metric on $(\Pi^*)^{-1}(\mathcal{P}_s)$. Thanks to this fact and the compactness of $(\Pi^*)^{-1}(V)$, the existence of such a $\zeta$  would imply that there exists $\zeta'(\epsilon,d,d^*,V)>0$ with $\zeta'\underset{\epsilon\rightarrow 0}{\longrightarrow} 0$ such that $$\sup_{n\in \mathbb{Z}\cap J_{exit}(x)}\big(d^*((\Phi^*_{pk})^n(x^*),(\Phi^*_{pk})^{n}(y^*))\big)< \zeta'$$ where  $x^*:=(\Pi^{*})^{-1}(x)$, $y^*:=(\Pi^*)^{-1}(y)$. Finally, thanks  to the continuity of $\Phi^*_{pk}$ and the compactness of $(\Pi^*)^{-1}(V)$, the above would imply that there exists $\widetilde{\eta}(\epsilon,d,d^*,V)>0$ with $\widetilde{\eta}\underset{\epsilon\rightarrow 0}{\longrightarrow} 0$ such that $$\sup_{t\in J_{exit}(x)}\big(d^*((\Phi_{pk}^*)^t(x^*),(\Phi_{pk}^*)^{t}(y^*))\big)< \widetilde{\eta}$$ We conclude the proof of Lemma \ref{l.comparisonmetricsflows} thanks to Claim 0. 
 \end{proof}
Fix for the rest of this proof a sufficiently small $\epsilon>0$, $x\in  V\cap (\mathcal{P}_0\setminus 0_{\mathcal{P}_0})$ and $y\in \mathcal{P}_0\setminus 0_{\mathcal{P}_0}$ such that 
\begin{equation}\label{eq.closexy}
    \sup_{t\in J_{exit}(x)}\big(d(\Phi_{pk}^t(x),\Phi_{pk}^{t}(y))\big)<\epsilon
\end{equation}

Recall that the first return map of $\Phi_{pk}$ on $\mathcal{P}_0$, namely $\Phi_{pk}^1$, coincides with the map $\phi_{pk}$ and the intersections of $F^s_{pk},F^u_{pk}$, the stable and unstable foliations of $\Phi_{pk}$, with $\mathcal{P}_0$ coincide with $\mathcal{F}_p^s,\mathcal{F}_p^u$, where $0_{\mathcal{P}_0}$ is the unique singular point of $\mathcal{F}_p^s,\mathcal{F}_p^u$. We would like to show that there exists $\zeta(\epsilon,d,d^*,V)>0$ with $\zeta\underset{\epsilon\rightarrow 0}{\longrightarrow} 0$ such that 
\begin{equation}\label{eq.goal}
    \sup_{n\in \mathbb{Z}\cap J_{exit}(x)}\big(d^{pol}_p(\phi_{pk}^n(x),\phi_{pk}^{n}(y))\big)< \zeta
\end{equation}

Notice first that thanks to Inequality \ref{eq.closexy}, if $x\in \mathcal{F}^s_p(0_{\mathcal{P}_0})$ (resp. $x\in \mathcal{F}^u_p(0_{\mathcal{P}_0})$), then $y\in \mathcal{F}^s_p(0_{\mathcal{P}_0})$ (resp. $y\in \mathcal{F}^u_p(0_{\mathcal{P}_0})$). Let $L_1^s,...,L_p^s,L_1^u,...,L_p^u$ be the set of stable and unstable prongs of $0_{\mathcal{P}_0}$ in $\mathcal{P}_0$. In the following lines, we will prove among others that if $x\in L_i^{s}$ (resp. $x\in L_i^{u}$), then we also have that $y\in L^{s}_i$ (resp. $y\in L_i^{u}$). Denote by $Q_1,...,Q_{2p}$ the quadrants of $0_{\mathcal{P}_0}$ inside $(\mathcal{P}_0,\mathcal{F}_p^s,\mathcal{F}_p^u)$.

\textbf{Claim $\text{N}^{o}$ 2:} There exists $\epsilon_0(d,V)$ such that if $\epsilon<\epsilon_0$, then for every $n\in \mathbb{Z}\cap J_{exit}(x)$ if $\phi_{pk}^{n}(x)\in Q_l$, where $l\in \llbracket 1, 2p \rrbracket$, then $\phi_{pk}^{n}(y)\in Q_l$. 
\begin{proof}
Suppose that $\phi_{pk}^{n}(x)$ and $\phi_{pk}^{n}(y)$ do not belong in the same quadrant of $0_{\mathcal{P}_0}$ for some $n\in \mathbb{Z}\cap J_{exit}(x)$. This implies that $x,y$ do not belong in the same quadrant of $0_{\mathcal{P}_0}$. By hypothesis, since $x,y$ are $\epsilon$-close, $x\in V$ and $V$ is compact, for $\epsilon$ sufficiently small, we have the three following possibilities: 
\begin{enumerate}
    \item $x,y$ are close to $0_{\mathcal{P}_0}$
    \item $x,y$ are close to the common stable boundary of two distinct quadrants of $0_{\mathcal{P}_0}$ and $x,y\notin \mathcal{F}^{s,u}_p(0_{\mathcal{P}_0})$
    \item $x,y$ are close to the common unstable boundary of two distinct quadrants of $0_{\mathcal{P}_0}$ and $x,y\notin \mathcal{F}^{s,u}_p(0_{\mathcal{P}_0})$
\end{enumerate}
Assume without any loss of generality that we are in the second case (the other cases are treated similarly). If $\epsilon$ is sufficiently small, since $V\cap \mathcal{P}_0$ is a standard transverse polygon, the positive orbit of $x$ by $\phi_{pk}$ will remain inside $V\cap \mathcal{P}_0$ for a long time. Take $N$ to be the biggest element inside $J_{exit}(x)\cap \mathbb{N}$. Since $0_{\mathcal{P}_0}$ is a $p$-prong singularity with $p\geq 2$, $\phi_{pk}^N(x)$ and $\phi_{pk}^N(y)$ are going to be close to two distinct unstable prongs of $0_{\mathcal{P}_0}$ (see Figure \ref{f.dynamicspolygon}) and furthermore, by definition of $N$, $\phi_{pk}^{N+1}(x)\notin V$. It is easy to see that there exists $\epsilon_0(d,V)$ such that if $\epsilon<\epsilon_0$, then the previous two conditions imply that $d(\phi_{pk}^N(x),\phi_{pk}^N(y))>\epsilon$, which contradicts Inequality \ref{eq.closexy}.

\begin{figure}[h!]
    \centering
    \includegraphics[scale=0.1]{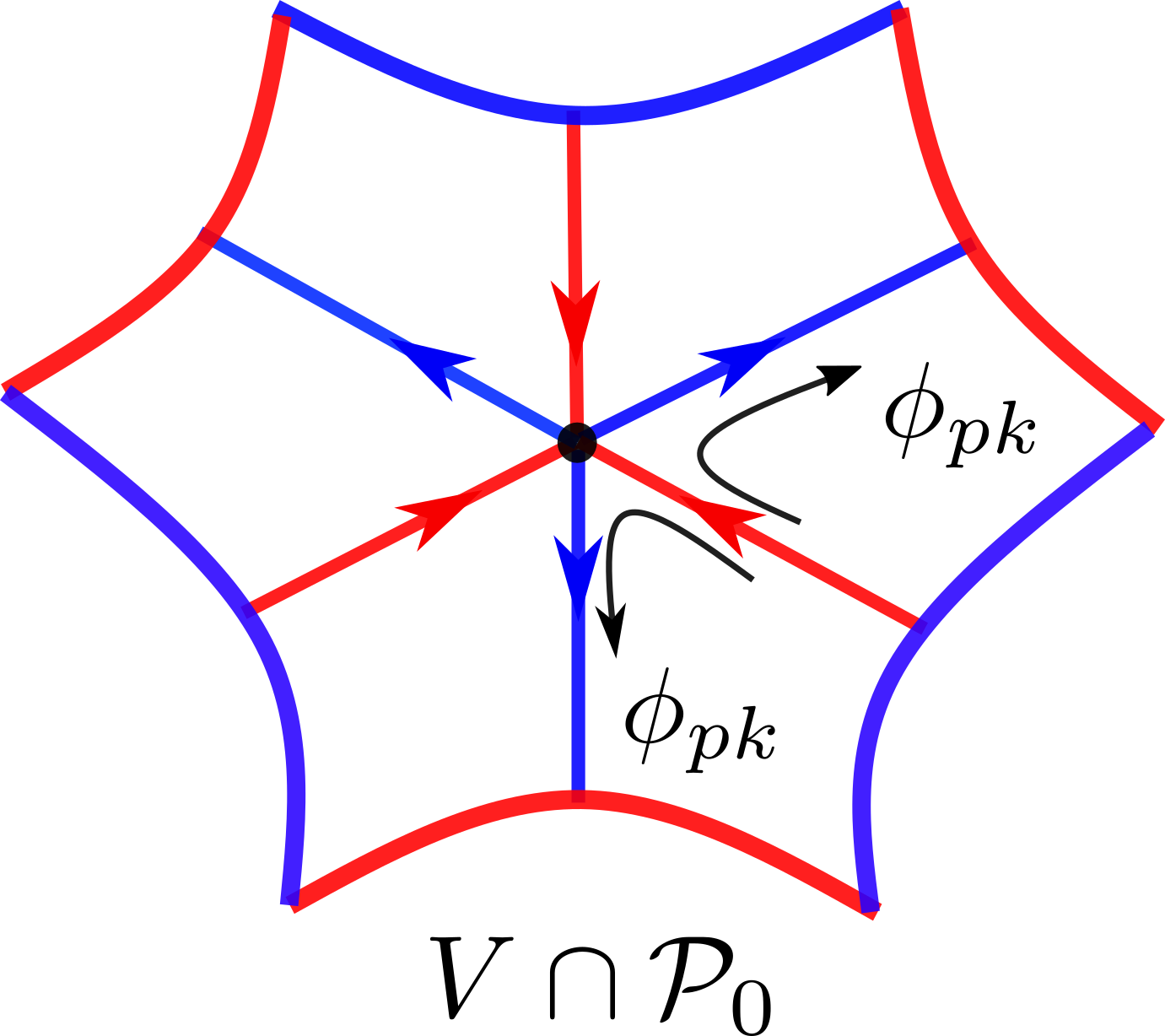}
    \caption{An illustration of the action of $\phi_{pk}$ when $p=3$ and $k=0$.}
    \label{f.dynamicspolygon}
\end{figure}. 
\end{proof}

If $\epsilon$ is sufficiently small and $x\in L^s_i$, then thanks to Claim 2, we also have that $y\in L^s_i$. In that case, using Items (1) and (5) of Proposition \ref{p.metricscomparison} and the facts that $x\in V$ and $V$ is compact, we get that  Inequality \ref{eq.closexy} gives \ref{eq.goal}, which finishes the proof of the lemma, thanks to Claim 1. Same result when $x\in L^u_i$. We can therefore suppose from now on that 
\begin{equation}\label{eq.hypothesis}
x,y\notin \mathcal{F}^{s,u}(0_{\mathcal{P}_0})\text{ and thus }J_{exit}(x)\text{ is compact.}
\end{equation}

Define $\pi_s:\mathcal{P}_0\rightarrow \mathcal{F}^s_p(0_{\mathcal{P}_0})$, $\pi_u:\mathcal{P}_0\rightarrow \mathcal{F}^s_p(0_{\mathcal{P}_0})$ as the unique continuous functions for which for every $z\in\mathcal{P}_0$ we have $\mathcal{F}_p^s(\pi_u(z))\cap \mathcal{F}_p^u(\pi_s(z))=\{z\}$. 

Also, for any $z,z'\in \mathcal{F}^s_p(0_{\mathcal{P}_0})$ denote by $[z,z']^s$ the unique stable segment in $\mathcal{F}^s_p(0_{\mathcal{P}_0})$ going from $z$ to $z'$. We similarly define $[z,z']^u$ for any two points $z,z'\in \mathcal{F}^u_p(0_{\mathcal{P}_0})$. We will say that a stable (resp. unstable) prong $L^s_i$ (resp. $L^u_i$) of $0_{\mathcal{P}_0}$ in $\mathcal{P}_0$ has \emph{period} $N$, if $N$ is the smallest positive integer for which $\phi_{pk}^{N}(L^s_i)\subset L^s_i$ (resp. $\phi_{pk}^{-N}(L^u_i)\subset L^u_i$). If $N$ is the period of $L^s_i$ (resp. $L^u_i$), then for any $z$ in $L^s_i\setminus 0_{\mathcal{P}_0}$ (resp. $L^u_i\setminus 0_{\mathcal{P}_0}$), we will call the segment $[z,\phi_{pk}^{N}(z)]^s$ (resp. $[z,\phi_{pk}^{-N}(z)]^u$), a \emph{prong fundamental domain}. 

\textbf{Claim $\text{N}^{o}$ 3:} There exists $\epsilon_1(d,V)$ such that if $\epsilon<\epsilon_1$, then for every $n\in \mathbb{Z}\cap J_{exit}(x)$ the segment $[\pi^s(\phi_{pk}^{n}(x)), \pi^s(\phi_{pk}^{n}(y))]^s$ (resp. $[\pi^u(\phi_{pk}^{n}(x)), \pi^u(\phi_{pk}^{n}(y))]^s$) does not contain a prong fundamental domain. 
\begin{proof}
    Assume that there exists $n\in \mathbb{Z}\cap J_{exit}(x)$ such that $[\pi^s(\phi_{pk}^{n}(x)), \pi^s(\phi_{pk}^{n}(y))]^s$ contains a prong fundamental domain. Consider $N\in -\mathbb{N}$ the smallest element in $\mathbb{Z}\cap J_{exit}(x)$ (this is well defined thanks to \ref{eq.hypothesis}). By definition of $N$ and since $\pi^s$ commutes with $\phi_{pk}$, we have that $[\pi^s(\phi_{pk}^{N}(x)), \pi^s(\phi_{pk}^{N}(y))]^s$ contains a prong fundamental domain and that $\phi_{pk}^{N-1}(x)\notin V$. Thanks to the two previous facts and the compactness of $V$, there exists $\epsilon_1(d,V)$ such that if $\epsilon<\epsilon_1$, then the two previous properties imply that $d(\phi_{pk}^{N}(x),\phi_{pk}^{N}(y))>\epsilon$, which contradicts Inequality \ref{eq.closexy} and finishes the proof of the claim. 
\end{proof}
Consider for every $i\in \llbracket 1, p\rrbracket$ two closed segments $ I_i^s\subset (L_i^s\setminus{0_{\mathcal{P}_0}})\cap V$ and $I_i^u\subset (L_i^u\setminus{0_{\mathcal{P}_0}})\cap V$ containing each at least two disjoint prong fundamental domains. Fix $U\subset V\cap \mathcal{P}_0$ a very small open neighborhood of $0_{\mathcal{P}_0}$ such that:
\begin{enumerate}
    \item $U$ is disjoint from $(\pi^s)^{-1}(I_i^s)$ and $(\pi^u)^{-1}(I_i^u)$ for every $i\in \llbracket 1, p\rrbracket$
    \item for every $z\in U\setminus \mathcal{F}^s_p(0_{\mathcal{P}_0})$, there exist $N\in\mathbb{N}^*$ and $i\in \llbracket 1,p\rrbracket$ such that for every $l\in \llbracket 0, N\rrbracket$ we have $\phi_{pk}^l(z)\in V\cap \mathcal{P}_0$ and $\pi^u(\phi_{pk}^N(z))\in I_i^u$
    \item for every $z\in U\setminus \mathcal{F}^u_p(0_{\mathcal{P}_0})$ there exist $M\in\mathbb{N}^*$ and $j\in \llbracket 1,p\rrbracket$ such that for every $l\in \llbracket 0, M\rrbracket$ we have $\phi_{pk}^{-l}(z)\in V\cap \mathcal{P}_0$ and  $\pi^s(\phi_{pk}^{-M}(z))\in I_j^s$
\end{enumerate}
The existence of such a $U$ follows from the fact that $V\cap\mathcal{P}_0$ is a standard transverse polygon for $\Phi_{pk}$. Fix $U'$ an open neighborhood of $0_{\mathcal{P}_0}$ such that $\text{Clos}(U')\subset U$ and take $\epsilon$ sufficiently small so that the $\epsilon$-neighborhood of $U'$ is contained in $U$. Moreover, by endowing $V-U'$ with the metric $d^{pol}_p$, thanks to the compactness of $V-U'$, there exists $\zeta_1(\epsilon,d, V)>0$ with  $\zeta_1\underset{\epsilon\rightarrow 0}{\longrightarrow}0$ such that for any $z\in V-U'$ and $ z'\in \mathcal{P}_0$ $$d(z,z')<\epsilon\implies d^{pol}_p(z,z')<\zeta_1$$
    
It follows that, for any $n\in \mathbb{Z}\cap J_{exit}(x)$ for which $\phi^n_{pk}(x)\notin U'$, we have that 
\begin{equation}\label{eq.almostresult}
    d^{pol}_p(\phi^n_{pk}(x),\phi^n_{pk}(y))<\zeta_1(\epsilon,d,V)
\end{equation}
where $\zeta_1\underset{\epsilon\rightarrow 0}{\longrightarrow}0$.

If for every $n\in \mathbb{Z}\cap J_{exit}(x)$ we have that $\phi^n_{pk}(x)\notin U'$, then \ref{eq.almostresult} finishes the proof of the lemma, thanks to Claim 1. Consider now the case where there exists $n\in \mathbb{Z}\cap J_{exit}(x)$ such that $\phi^n_{pk}(x)\in U'$. By eventually replacing $x$ by  $\phi^n_{pk}(x)$, assume without any loss of generality that $x\in U'$. By the definition of $U'$ and $U$, we have that $x,y\in U$ and thanks to Claims 2 and 3 we also have that: 
\begin{enumerate}
    \item if $\epsilon$ is sufficiently small, $x,y$ belong in the same quadrant of $0_{\mathcal{P}_0}$ inside  $(\mathcal{P}_0,\mathcal{F}_p^s,\mathcal{F}_p^u)$
    \item there exist $N\in\mathbb{N}^*$ and $i\in \llbracket 1,p\rrbracket$ such that for every $l\in \llbracket 0, N\rrbracket$ we have $\phi_{pk}^l(x),\phi_{pk}^l(y)\in V\cap \mathcal{P}_0$ and $\pi^u(\phi_{pk}^N(x)), \pi^u(\phi_{pk}^N(y))\in I_i^u$
    \item there exist $M\in\mathbb{N}^*$ and $j\in \llbracket 1,p\rrbracket$ such that for every $l\in \llbracket 0, M\rrbracket$ we have $\phi_{pk}^{-l}(x),\phi_{pk}^{-l}(y)\in V\cap \mathcal{P}_0$ and $\pi^s(\phi_{pk}^{-M}(x)), \pi^s(\phi_{pk}^{-M}(y))\in I_j^s$
\end{enumerate}

By using the Inequalities \ref{eq.closexy},  \ref{eq.almostresult} and Proposition \ref{p.metricscomp2}, we get that if $\epsilon$ is sufficiently small, there exists $\zeta(\epsilon,d,V)$ with $\zeta\underset{\epsilon\rightarrow 0}{\longrightarrow}0$, such that for every $n\in \mathbb{Z}\cap J_{exit}(x)$ we have that $d^{pol}_p(\phi^n_{pk}(x),\phi^n_{pk}(y))<\zeta(\epsilon,d,V)$. This finishes the proof of the lemma, thanks to Claim 1. 
\end{proof}
\subsection{Proof of Proposition \ref{p.stablestar}}\label{s.proofproposition}
Consider $M$ a closed 3-manifold, $\Phi=(\Phi^t)_{t\in \mathbb{R}}$ an almost pseudo-Anosov flow on $M$, $F^{s,u}$ the stable and unstable foliations of $\Phi$ and $\gamma$ a periodic orbit of $\Phi$, whose tubular neighborhoods are solid tori. Let $U_\gamma$ be the (solid torus) neighborhood of $\gamma$, $(N_{pk},\Phi_{pk})$ be the hyperbolic model and $H:U_\gamma\rightarrow V_{pk}\subset N_{pk}$ the orbital equivalence associated to $\gamma$ by Proposition \ref{p.aroundcircleprong}, where  $p\geq 2$, $k\in \llbracket 0, p-1\rrbracket$ and $V_{pk}$ is a neighborhood of $\gamma_{pk}:=H(\gamma)$. Recall that $$N_{pk}:= \frac{\mathbb{R}^2\times \mathbb{R}}{((x,y),t+1) \sim (\phi_{pk}(x,y),t) }$$ 
For any $t\in\mathbb{R}$, denote by $\mathcal{P}_t$ the projection of $\mathbb{R}^2\times \{t\}$ on $N_{pk}$ and by $0_{\mathcal{P}_s}$ the unique point in $\mathcal{P}_s\cap \gamma_{pk}$. Consider $(N^*_{pk},\Phi^*_{pk})$ the blow-up of $(N_{pk},\Phi_{pk})$ along $\gamma_{pk}$, $\Pi^*$ its associated blow-down map and $(F^{s,u}_{pk})$ its associated stable and unstable foliations. Consider also $(M^*, \Phi^*)$ a flow obtained after blowing up $\gamma$, $\Pi^*_M$ its associated  blow-down map and $(F^{s,u})^*$ its associated stable and unstable foliations. Recall that by our construction of $(M^*,\Phi^*)$ in Section \ref{s.blowupflows}, 
\begin{rema}\label{r.recall}\quad

\begin{itemize}
    \item $M^*$ was obtained by removing $U_\gamma$ from $M$ and glueing at its place $V^*_{pk}$
    \item $(\Pi_M^*)^{-1}(U_\gamma)=V_{pk}^*$ 
    \item the restriction of $\Pi_M^*$ on $(\Pi_M^*)^{-1}(U_\gamma)$ coincides with $H^{-1}\circ \Pi^*: V_{pk}^*\rightarrow U_\gamma$ 
    \item $\Phi^*$ inside $(\Pi_M^*)^{-1}(U_\gamma)$ is obtained by reparametrizing $\Phi_{pk}^*$ inside $V_{pk}^*$
    \item $\Pi^*_M$ defines an orbital equivalence between $(\Phi^*, \inte{M^*})$ and $(\Phi, M-\gamma)$
\end{itemize} 
\end{rema}
This section will be dedicated to the proof of Proposition \ref{p.stablestar}:
\begin{prop*}For any metric $d^*$ on $M^*$ (compatible with the the topology of $M^*$), any $x^*\in M^*$ contained in a stable (resp. unstable) leaf in $(F^s)^*$ (resp. $(F^u)^*$) and any $y^*\in (F^s)^*(x)$ (resp. $y^*\in (F^u)^*(x)$), there exists an increasing homeomorphism $h:\mathbb{R}\rightarrow \mathbb{R}$ such that 
\begin{equation*} d^*((\Phi^*)^t(x^*),(\Phi^*)^{h(t)}(y^*))\underset{t\rightarrow +\infty}{\longrightarrow} 0~~ \big(\text{resp. }d^*((\Phi^*)^t(x^*),(\Phi^*)^{h(t)}(y^*))\underset{t\rightarrow -\infty}{\longrightarrow}0\big)
\end{equation*}
\end{prop*}
Before beginning the proof of the above proposition, let us first state the following useful variant of Lemma \ref{l.comparisonmetricsflows}:

\begin{lemm}\label{l.keylemma}
    Take $U\subset U_\gamma$ a neighborhood of $\gamma$ homeomorphic to a solid torus such that for every $s\in [0,1]$ the set $H(U)\cap \mathcal{P}_s$ is a transverse standard polygon in $N_{pk}$. For every $z\in U$, let $T^-(z):=\inf\{t<0|\Phi^{t}(x)\in U\}$, $T^+(z):=\sup\{t>0|\Phi^{t}(x)\in U\}$ and $J_{exit}(z)$ be the closure in $\mathbb{R}$ of the interval $(T^-(z), T^+(z))$.

    For any $\epsilon>0$ sufficiently small, there exist $\eta''(\epsilon,d,d^*,U)>0$ with $\eta'' \underset{\epsilon\rightarrow 0}{\longrightarrow} 0$ such that for any $z\in U-\gamma$, any continuous and strictly increasing map $h:J_{exit}(z)\rightarrow \mathbb{R}$ and any $w\in M-\gamma$,  $$\sup_{t\in J_{exit}(z)}\big(d(\Phi^t(z),\Phi^{h(t)}(w))\big)< \epsilon \implies  \sup_{t\in J_{exit}(z)}\big(d^*((\Pi_M^*)^{-1}(\Phi^t(z)),(\Pi_M^*)^{-1}(\Phi^{h(t)}(w)))\big)< \eta''$$
\end{lemm}
The above lemma is an immediate consequence of Lemma \ref{l.comparisonmetricsflows} and Remark \ref{r.recall}. We will thus omit its proof.

\begin{proof}[Proof of Proposition \ref{p.stablestar}]
Fix  $d$ a metric on $M$, $d^*$ a metric on $M^*$, $x^*\in M^*$ contained in a stable leaf in $(F^s)^*$ and $y^*\in (F^s)^*(x)$. Let $x:=\Pi^*_M(x^*)$ and $y:=\Pi^*_M(y^*)\in F^s(x)$. 

Consider first the case where $x,y \in F^s(\gamma)$. In this case, there exists $\mathcal{O}^*$ a periodic orbit of $\Phi^*$ on $\partial M^*$ such that $x^*,y^*\in (F^s)^*(\mathcal{O}^*)$. Since $\Pi^*_M$ defines an orbital equivalence between $(\Phi^*, \inte{M^*})$ and $(\Phi, M-\gamma)$, by replacing if necessary $x^*,y^*$ by $(\Phi^*)^T(x^*),(\Phi^*)^{T'}(y^*)$ for some large $T,T'>0$, we can assume that $x,y\in U_\gamma$ and that the positive orbits of $x,y$ by $\Phi$ remain in $U_\gamma$. This implies that the positive orbits of both $x^*$ and $y^*$ remain inside $(\Pi_M^*)^{-1}(U_\gamma)$ and thanks to Remark \ref{r.recall} the previous orbits both accumulate to $\mathcal{O}^*$, which trivially gives the desired result. 

Assume from now on that $x,y\notin F^s(\gamma)$. By Definition \ref{d.pseudoanosovflow}, there exists an increasing homeomorphism $h_M:\mathbb{R}\rightarrow \mathbb{R}$ such that 
\begin{equation}\label{eq.convergencehypothesis}
    d(\Phi^t(x),\Phi^{h_M(t)}(y))\underset{t\rightarrow +\infty}{\longrightarrow} 0
\end{equation} Fix $\epsilon>0$. By replacing if necessary $x,y$ by $\Phi^T(x),\Phi^{h(T)}(y)$ for some $T>0$, we can assume without any loss of generality that $d(\Phi^t(x),\Phi^{h_M(t)}(y))<\epsilon $ for every $t>0$. 

Take $U\subset U_\gamma$ a neighborhood of $\gamma$ homeomorphic to a solid torus such that $x\notin U$ and such that for every $s\in [0,1]$ the set $H(U)\cap \mathcal{P}_s$ is a transverse standard polygon in $N_{pk}$. Notice that thanks to the continuity of $\Pi_M^*$, there exists $\eta(\epsilon,d,d^*,U)>0$ with $\eta\underset{\epsilon\rightarrow 0}{\longrightarrow}0$ such that for every $z,z'\in \text{Clos}(M-U)$ 
\begin{equation}\label{eq.bound0}
   d(z,z')<\epsilon \implies d^*((\Pi_M^*)^{-1}(z),(\Pi_M^*)^{-1}(z'))<\eta 
\end{equation}

Next, consider the sequence of intervals $(I_i)_{i\in\mathbb{N}}=([a_i,b_i])_{i\in\mathbb{N}}$ defined by the following properties (recall that $x,y\notin F^s(\gamma)$): 
\begin{itemize}
    \item $a_0=0$, $b_i=a_{i+1}$ for every $i\in \mathbb{N}$ 
    \item for every $i\in 2\mathbb{N}$ and every $t\in (a_i,b_i)$ we have that $\Phi^t(x)\notin U$ 
    \item for every $i\in 2\mathbb{N}+1$ and every $t\in [a_i,b_i]$ we have that $\Phi^t(x)\in U$ 
\end{itemize}
Thanks to \ref{eq.convergencehypothesis} and \ref{eq.bound0}, there exists a positive sequence $(\eta_i)_{i\in 2\mathbb{N}}$ with $\eta_i\underset{i\rightarrow +\infty}{\longrightarrow}0$ such that for every $i\in 2\mathbb{N}$ and every $t\in I_i$ 
\begin{equation}\label{eq.step1}
    d^*((\Pi_M^*)^{-1}(\Phi^t(x)),(\Pi_M^*)^{-1}(\Phi^{h_M(t)}(y)))<\eta_i
\end{equation}
Thanks to \ref{eq.convergencehypothesis} and Lemma \ref{l.keylemma}, if $\epsilon$ is sufficiently small, there exists a positive sequence $(\eta''_i)_{i\in 2\mathbb{N}+1}$ with $\eta''_i\underset{i\rightarrow +\infty}{\longrightarrow}0$ such that for every $i\in 2\mathbb{N}+1$ and for every $t\in I_i$ 
\begin{equation}\label{eq.step2}
      d^*((\Pi_M^*)^{-1}(\Phi^t(x)),(\Pi_M^*)^{-1}(\Phi^{h_M(t)}(y)))< \eta_i''
\end{equation}
 By combining \ref{eq.step1} with \ref{eq.step2} and using the fact that $\Pi^*_M$ defines an orbital equivalence between $(\Phi^*, \inte{M^*})$ and $(\Phi, M-\gamma)$, there exists a continuous and strictly increasing map $h:\mathbb{R}\rightarrow \mathbb{R}$ such that $d^*((\Phi^*)^{t}(x^*),(\Phi^*)^{h(t)}(y^*))\underset{t\rightarrow +\infty}{\longrightarrow}0$, which finishes the proof of the proposition.  
\end{proof}
\subsection{Proof of Theorem \ref{t.maintheorem}}
Consider $M$ a closed 3-manifold, $d$ a metric on $M$, $\Phi=(\Phi^t)_{t\in \mathbb{R}}$ an expansive flow on $M$, $F^{s,u}$ its stable and unstable foliations  and $\gamma$ a periodic orbit of $\Phi$, whose tubular neighborhoods are  solid tori. Let ($M^*$,$ \Phi^*$) be a flow obtained after blowing up $\gamma$, $\Pi^*_M$ its associated blow-down map and $(F^{s,u})^*$ its associated stable and unstable foliations. Thanks to Lemma \ref{l.goodfoli}, by eventually reparametrizing $\Phi^*$, we get that $\Phi^*$ admits a good foliation $\mathcal{F}$ on $\partial M^*$. Denote by $(M_\mathcal{F},\Phi_\mathcal{F})$ the flow obtained by crushing every leaf of $\mathcal{F}$ to a point, by $\Pi_{\mathcal{F}}$ the natural projection from $M^*$ to $M_\mathcal{F}$ and by $F^{s,u}_\mathcal{F}$ the stable and unstable foliations given by Proposition-Definition \ref{p.foliationsaftersurgery}. Our goal in this section consists in proving that $\Phi_\mathcal{F}$ is expansive if and only if $F^{s,u}_\mathcal{F}$ does not admit an 1-prong singularity. 

The fact that if $\Phi_\mathcal{F}$ is expansive, then $F^{s,u}_{\mathcal{F}}$ does not admit an 1-prong singularity follows immediately from Theorem \ref{t.expansiveimpliespseudo}. Let us now prove the converse. Fix $d_\mathcal{F}$ a metric on $M_{\mathcal{F}}$ (compatible with its topology).

Suppose that there exist two positive sequences $\epsilon_m \longrightarrow 0$, $t_m\longrightarrow +\infty$, a sequence of increasing homeomorphisms $h_m:\mathbb{R}\rightarrow \mathbb{R}$ with $h_m(0)=0$ and two sequences of points $x_m, y_m\in M$ such that $y_m\notin \underset{t\in[-t_m,t_m]}{\cup}\Phi_{\mathcal{F}}^t(x_m)$ and $d_\mathcal{F}(\Phi_\mathcal{F}^t(x_m),\Phi_\mathcal{F}^{h_m(t)}(y_m))<\epsilon_m$ for all $t\in \mathbb{R}$ and $m\in \mathbb{N}$.

Recall that, by our discussion in Section \ref{s.dehn-frieddefi}, $\Phi_\mathcal{F}$ is an almost pseudo-Anosov flow. Using this fact and Proposition \ref{p.inversesurgery}, we can perform a Dehn-Fried surgery on $\gamma_\mathcal{F}$ in order to obtain a flow orbitally equivalent to $\Phi$. Let ($M_\mathcal{F}^*$,$ \Phi_\mathcal{F}^*$) be a flow obtained after blowing up $\gamma_\mathcal{F}$, $\Pi^*_\mathcal{F}$ its associated blow-down map, $(F_\mathcal{F}^{s,u})^*$ its associated stable and unstable foliations and $d_\mathcal{F}^*$ any metric on $M_\mathcal{F}^*$ (compatible with its topology).

Thanks to Proposition \ref{p.aroundcircleprongaftersurgery}, neither $x_m$ nor $y_m$ can belong to $\gamma_\mathcal{F}$ for $m$ sufficiently big. Take $m$ sufficiently big so that $x_m,y_m\notin \gamma_\mathcal{F}$. Let $U_{\gamma_\mathcal{F}}$ be the (solid torus) neighborhood of $\gamma$, $(N_{pk},\Phi_{pk})$ be the hyperbolic model and $H:U_\gamma\rightarrow V_{pk}\subset N_{pk}$ the orbital equivalence associated to $\gamma$ by Proposition \ref{p.aroundcircleprong}, where  $p\geq 2$, $k\in \llbracket 0, p-1\rrbracket$ and $V_{pk}$ is a neighborhood of $\gamma_{pk}:=H(\gamma_\mathcal{F})$. Consider for every $s\in \mathbb{R}$ the plane $\mathcal{P}_s\subset N_{pk}$ defined in the beginning of the previous section and take $U\subset U_{\gamma_\mathcal{F}}$ a  neighborhood of $\gamma_\mathcal{F}$ homeomorphic to a solid torus such that for every $s\in [0,1]$ the set $H(U)\cap \mathcal{P}_s$ is a transverse standard polygon in $N_{pk}$.

For any $\epsilon$ sufficiently small, thanks to the continuity of $\Pi_\mathcal{F}^*$, there exists $\eta(\epsilon,d_\mathcal{F},d_\mathcal{F}^*,U)>0$ with $\eta\underset{\epsilon\rightarrow 0}{\longrightarrow}0$ such that for every $z,z'\in \text{Clos}(M-U)$ 
\begin{equation}\label{eq.bound1}
   d_\mathcal{F}(z,z')<\epsilon \implies d_\mathcal{F}^*((\Pi_\mathcal{F}^*)^{-1}(z),(\Pi_\mathcal{F}^*)^{-1}(z'))<\eta 
\end{equation}
Assume for the sake of simplicity that $x_m\notin F^s_\mathcal{F}(\gamma_\mathcal{F})\cup F^u_\mathcal{F}(\gamma_\mathcal{F})$ (the case where $x_m\in F^s_\mathcal{F}(\gamma_\mathcal{F})\cup F^u_\mathcal{F}(\gamma_\mathcal{F})$ can be treated by a similar argument). Consider for every $m\in \mathbb{N}$ the sequence of intervals $(I^m_i)_{i\in\mathbb{Z}}=([a^m_i,b^m_i])_{i\in \mathbb{Z}}$ defined by the following properties: 
\begin{itemize}
    \item $b^m_i=a^m_{i+1}$ for every $i\in \mathbb{Z}$ 
    \item for every $i\in 2\mathbb{Z}$ and every $t\in (a^m_i,b^m_i)$ we have that $\Phi_\mathcal{F}^t(x_m)\notin U$ 
    \item for every $i\in 2\mathbb{Z}+1$ and every $t\in [a^m_i,b^m_i]$ we have that $\Phi_\mathcal{F}^t(x_m)\in U$ 
\end{itemize}
    Thanks to \ref{eq.bound1}, there exists a positive sequence $(\eta_m)_{m\in \mathbb{Z}}$ with $\eta_m\underset{m\rightarrow +\infty}{\longrightarrow}0$ such that for every $m\in \mathbb{N}$, every $i\in 2\mathbb{Z}$ and every $t\in I^m_i$ 
\begin{equation}\label{eq.step3}
    d_\mathcal{F}^*((\Pi_\mathcal{F}^*)^{-1}(\Phi_\mathcal{F}^t(x_m)),(\Pi_\mathcal{F}^*)^{-1}(\Phi_\mathcal{F}^{h_m(t)}(y_m)))<\eta_m
\end{equation}
By definition of $U$ and since $\Phi_\mathcal{F}$ is an almost pseudo-Anosov flow, we can apply Lemma \ref{l.keylemma} for $U$. Thanks to the previous lemma, there exists a positive sequence $(\eta''_m)_{m\in \mathbb{N}}$ with $\eta''_m\underset{m\rightarrow +\infty}{\longrightarrow}0$ such that for every  $m\in \mathbb{N}$ sufficiently big, every $i\in 2\mathbb{Z}+1$ and every $t\in I^m_i$ 
\begin{equation}\label{eq.step4}
      d_\mathcal{F}^*((\Pi_\mathcal{F}^*)^{-1}(\Phi_\mathcal{F}^t(x_m)),(\Pi_\mathcal{F}^*)^{-1}(\Phi_\mathcal{F}^{h_m(t)}(y_m)))<\eta_m''
\end{equation}
Let $x_m^*:=(\Pi_\mathcal{F}^*)^{-1}(x_m)$ and $y_m^*:=(\Pi_\mathcal{F}^*)^{-1}(y_m)$. By combining \ref{eq.step3} with \ref{eq.step4} and using the fact that $\Pi^*_\mathcal{F}$ defines an orbital equivalence between $(\Phi_\mathcal{F}^*, \inte{M_\mathcal{F}^*})$ and $(\Phi_\mathcal{F}, M_\mathcal{F}-\gamma_\mathcal{F})$,  we obtain that there exists two positive sequences $\eta^*_m \longrightarrow 0$, $t^*_m\longrightarrow +\infty$ and a sequence of increasing homeomorphisms $h^*_m:\mathbb{R}\rightarrow \mathbb{R}$ with $h^*_m(0)=0$ such that ${y^*_m}\notin \underset{t\in[-t^*_m,t^*_m]}{\cup}(\Phi^*_{\mathcal{F}})^t({x^*_m})$ and $d^*_\mathcal{F}((\Phi^*_{\mathcal{F}})^t({x^*_m}),(\Phi^*_{\mathcal{F}})^{h^*_m(t)}({y^*_m}))<\eta^*_m$ for all $t\in \mathbb{R}$ and $m\in \mathbb{N}$.

By reparametrizing $\Phi^*_{\mathcal{F}}$ if necessary, thanks to Proposition \ref{p.inversesurgery}, we have that $\Phi^*_\mathcal{F}$ admits a good foliation $\mathcal{G}$ such that if $(M_\mathcal{G},\Phi_\mathcal{G})$ is the flow obtained after crushing every leaf of $\mathcal{G}$ into a point, then $(M_\mathcal{G},\Phi_\mathcal{G})$ is orbitally equivalent to $(M,\Phi)$. Let $\Pi_\mathcal{G}$ be the projection from $M_{\mathcal{F}}^*$ to $M_{\mathcal{G}}\approx M$, $X_m:=\Pi_\mathcal{G}(x^*_m)$ and $Y_m:=\Pi_\mathcal{G}(y^*_m)$. Since  $\Pi_\mathcal{G}$ is continuous and $M$ and $M^*_{\mathcal{F}}$ are compact, we have that there exists two positive sequences $\zeta_m \longrightarrow 0$,  $s_m\longrightarrow +\infty$ and a sequence of increasing homeomorphisms $f_m:\mathbb{R}\rightarrow \mathbb{R}$ with $f_m(0)=0$ such that $Y_m\notin \underset{t\in[-s_m,s_m]}{\cup}\Phi^t(X_m)$ and $d(\Phi^t(X_m),\Phi^{f_m(t)}(Y_m))<\zeta_m$ for all $t\in \mathbb{R}$ and $m\in \mathbb{N}$. This contradicts the fact that $\Phi$ is an expansive flow, which finishes the proof of Theorem \ref{t.maintheorem}.

\end{document}